\documentclass[onecolumn]{autart}    

\usepackage{graphicx}
\usepackage{subfig}
\usepackage{amsmath}
\usepackage{amssymb}
\usepackage[font=footnotesize]{caption} 
\usepackage{subfig} 
\usepackage[noadjust]{cite} 
\usepackage{color}
\usepackage{algorithm} 
\usepackage{algpseudocode} 
\usepackage{enumerate}
\usepackage{setspace}

\usepackage{tikz}
\usetikzlibrary{calc} 
\usetikzlibrary{shapes} 
\usetikzlibrary{chains}
\usetikzlibrary{fit}
\usetikzlibrary{arrows}
\usetikzlibrary{decorations.text} 
\usetikzlibrary{decorations.markings}

\newtheorem{theorem}{Theorem}
\newtheorem{lemma}{Lemma}

\newtheorem{remark}{Remark}

\newtheorem{corollary}{Corollary}

\newtheorem{definition}{Definition}
\newtheorem{problem}{Problem}

\newcommand{\Null}{\mathrm{Null}}
\newcommand{\Range}{\mathrm{Range}}
\newcommand{\one}{\mathbf{1}}
\newcommand{\rank}{\mathrm{Rank}}
\newcommand{\mydiag}{\mathrm{diag}}

\newcommand{\blkdiag}{\mathrm{blkdiag}}

\newcommand{\T}{\mathrm{T}}

\newcommand{\R}{\mathbb{R}}
\newcommand{\G}{\mathcal{G}}
\newcommand{\E}{\mathcal{E}}
\newcommand{\V}{\mathcal{V}}
\newcommand{\N}{\mathcal{N}}

\newcommand{\B}{\mathcal{B}}

\newcommand{\dia}[1]{\mathrm{diag}\left(#1\right)} 
\newcommand{\myspan}[1]{\mathrm{span}\left\{ #1 \right\}}

\graphicspath{{figures/}}

\usepackage{soul}

\begin{document}

\begin{frontmatter}

\title{Localizability and Distributed Protocols for Bearing-Based Network Localization in Arbitrary Dimensions}

\author[shiyu]{Shiyu Zhao}\ead{shiyuzhao@engr.ucr.edu},    
\author[dan]{Daniel Zelazo}\ead{dzelazo@technion.ac.il}              
\address[shiyu]{Department of Mechanical Engineering, University of California, Riverside, USA.}  \address[dan]{Faculty of Aerospace Engineering, Technion - Israel Institute of Technology, Haifa, Israel.}

\begin{abstract}                          
This paper addresses the problem of bearing-based network localization, which aims to localize all the nodes in a static network given the locations of a subset of nodes termed anchors and inter-node bearings measured in a common reference frame.
The contributions of the paper are twofold.
Firstly, we propose necessary and sufficient conditions for network localizability with both algebraic and rigidity theoretic interpretations.
The analysis of the localizability heavily relies on the recently developed bearing rigidity theory and a special matrix termed the {bearing Laplacian}.
Secondly, we propose a linear distributed protocol for bearing-based network localization.
The protocol can globally localize a network if and only if the network is localizable.
The sensitivity of the protocol to constant measurement errors is also analyzed.
One novelty of this work is that the localizability analysis and localization protocol are applicable to networks in arbitrary dimensional spaces.
\end{abstract}
\begin{keyword}                           
Sensor network, Network localizability, Distributed localization, Bearing rigidity, Bearing Laplacian.
\end{keyword}                             

\end{frontmatter}

\section{Introduction}
Distributed localization of sensor networks is a core problem in many multi-agent coordination tasks.
\emph{Network localizability} and \emph{distributed protocols} are two fundamental problems for any network localization problems.
Network localizability characterizes whether or not a network can be possibly localized given the anchor locations and inter-neighbor relative measurements, whereas distributed protocols are used for localizing the network in a distributed manner if the network is localizable.
According to the types of the relative measurements used for localization, the existing works can be divided into three classes: distance-based, bearing-based, and position-based.
Distance-based network localization has been studied extensively so far (see \cite{AspnesTMC2006,Mao2007NetworkLocalization,Khan2009TSP,LinZhiyun2014TSP} and the references therein).
The analysis of the localizability in distance-based network localization relies heavily on the distance rigidity theory. It has been shown that a network in an $n$-dimensional space can be uniquely localized if the network is globally rigid and has at least $n+1$ anchors in a general position \cite{AspnesTMC2006}.
More recently, bearing-based network localization has also attracted extensive research attention \cite{BishopTAES2009,Niculescu2003AOA,ErenTurkishJournal2007,Piovan2013Automatica,Shames2013TAC,ZhuGuangwei2014Automatica,LinZhiyun2014ICCA}.
The analysis of the localizability in bearing-based network localization relies on the analogous bearing rigidity theory \cite{bishopconf2011rigid,Eren2012IJC,zelazo2014SE2Rigidity,zhao2014TACBearing}.
Finally, position-based network localization, where the inter-neighbor distance and local bearing measurements are used together for network localization, has been studied in \cite{LinZhiyun2015TSP} by using a complex graph Laplacian.

Although bearing-based network localization has been studied by many researchers, the two fundamental problems, network localizability and distributed protocols, have not yet been fully explored.
It was shown in \cite{ErenTurkishJournal2007,Piovan2013Automatica,Shames2013TAC,ZhuGuangwei2014Automatica} that a network is localizable when the network is bearing rigid and has at least two anchors.
This condition is, however, sufficient but not necessary when the number of anchors is greater than two \cite[Cor~10]{ZhuGuangwei2014Automatica}.
A necessary and sufficient condition for network localizability was proposed in \cite[Thm~15]{ZhuGuangwei2014Automatica} based on the notion of a stiffness matrix.
This condition is, however, applicable only to networks in two-dimensional spaces.
In fact, the localizability of a network is jointly determined by many factors such as its topological and Euclidean structure, as well as the selection of the anchors.
The relationship between the localizability and these factors have not been fully understood yet up to now.
Moreover, the existing bearing-based localization protocols are mainly applicable to networks in two-dimensional ambient spaces \cite{ErenTurkishJournal2007,Piovan2013Automatica,Shames2013TAC,ZhuGuangwei2014Automatica}.
General results of localizability or distributed protocols for bearing-based network localization in three and higher dimensional spaces are still lacking.

This paper studies the localizability and distributed protocols for bearing-based network localization in arbitrary dimensional spaces.
It is assumed that the anchors' locations and inter-neighbor bearings measured in a global reference frame are already given.
The main contributions of this work are summarized below.
\begin{enumerate}[(a)]
\item We first show that the bearing-based network localization problem can be formulated as a linear least-squares optimization problem.
    A special matrix termed the \emph{bearing Laplacian}, which can be viewed as a matrix-weighted graph Laplacian, emerges as a key part in the least-squares formulation and plays important roles in the subsequent analysis.
\item Based on the least-squares formulation, we propose necessary and sufficient conditions for network localizability with both algebraic and rigidity theoretic interpretations.
    These conditions not only provide numerical ways to examine the localizability of a given network but also provide intuitions on what a localizable network looks like.
\item We then propose a distributed linear localization protocol.
    It is proved that the protocol can globally localize a network if and only if the network is localizable.
    The sensitivity of the protocol to constant measurement errors is also analyzed.
\end{enumerate}
Finally, it is worth noting that the localizability analysis presented in this paper is independent to whether the sensing graph is directed or undirected. The convergence analysis of the proposed localization protocol, however, relies on the assumption of undirected sensing graphs. The convergence analysis of the protocol in the directed case is considered in \cite{zhao2015CDC}.

The rest of the paper is organized as follows.
Section~\ref{section_problemStatement} presents the linear least-squares formulation of the bearing-based network localization problem.
Section~\ref{section_propertiesOfBearingLaplacian} analyzes the properties of the bearing Laplacian and its connection to the bearing rigidity theory.
Section~\ref{section_localizabilityAnalysis} presents necessary and sufficient conditions for network localizability.
Section~\ref{section_distributedProtocol} proposes and analyzes a linear distributed localization protocol.
Conclusions are drawn in Section~\ref{section_conclusion}.

\paragraph*{Notations:} Given $A_i\in\mathbb{R}^{p\times q}$ for $i=1,\dots,n$, denote $ \mydiag(A_i)\triangleq\blkdiag\{A_1,\dots,A_n\}\in\mathbb{R}^{np\times nq}$.
Let $\|\cdot\|$ be the Euclidian norm of a vector or the spectral norm of a matrix, and $\otimes$ be the Kronecker product.
Denote $I_d\in\R^{d\times d}$ as the identity matrix, and $\one_d\triangleq[1,\dots,1]^\T\in\R^d$.
Let $\Null(\cdot)$ and $\Range(\cdot)$ be the null space and range space of a matrix, respectively.

\section{Problem Formulation of Bearing-Based Network Localization}
\label{section_problemStatement}

In this section, the problem of bearing-based network localization is formally stated and then formulated as a linear least-squares problem.
Central to this problem is the notion of localizability, which is formally defined here.

\subsection{Problem Statement}

Consider a network of $n$ stationary nodes in $\R^d$ ($n\ge2$ and $d\ge2$).
Assume no two nodes are collocated.
Let $p_i \in \R^d$ be the location of node $i$ ($i=1,\dots,n$).
Define the \emph{edge vector} and the \emph{bearing} between nodes $i$ and $j$ as
\begin{align*}
e_{ij}\triangleq p_j-p_i, \quad g_{ij}\triangleq \frac{e_{ij}}{\|e_{ij}\|}.
\end{align*}
The unit vector $g_{ij}$ represents the relative bearing of $p_j$ with respect to $p_i$.
Note $e_{ij}=-e_{ji}$ and $g_{ij}=-g_{ji}$.
Suppose the locations of $n_a$ \emph{anchor} nodes are already given and the locations of the remaining $n_f$ \emph{follower} nodes are to be estimated ($n_a+n_f=n$).
Denote $\V_a=\{1,\dots,n_a\}$, $\V_f=\{n_a+1,\dots,n\}$, and $\V=\V_a\cup\V_f$.
Denote $p_a=[p_1^\T ,\dots,p_{n_a}^\T]^\T\in\R^{dn_a}$, $p_f=[p_{n_a+1}^\T ,\dots,p_{n}^\T ]^\T\in\R^{dn_f}$, and $p=[p_a^\T ,p_f^\T]^\T\in\R^{dn}$.

Suppose each node has the bearing-only sensing capabilities.
The sensing topology of the network defines a graph $\G=(\V,\E)$ where $\E\subset\mathcal{V} \times \mathcal{V}$.
Denote $(i,j)$ as the directed edge with node $i$ as the tail and node $j$ as the head.
The directed edge $(i,j)\in\E$ indicates that node $i$ can ``see'' node $j$; that is node $i$ can measure the relative bearings $g_{ij}$ of node $j$.
Node $j$ is called the neighbor of node $i$ if $(i,j)\in\E$, and $\mathcal{N}_i\triangleq\{j \in \mathcal{V}|  (i,j)\in \mathcal{E}\}$ is the neighbourhood of node $i$.
We assume a {global orientation} that can be sensed by all the nodes, and thus all measured bearings can be expressed with respect to this common orientation.
The global orientation means a common north for the two-dimensional space, and a common north-east-down reference for the three-dimensional space.
Finally, let $\G(p)$ denote the network that is the graph $\G$ with each vertex $i \in \V$ mapped to the point $p_i$.

The problem of bearing-based network localization is formally stated below.

\begin{problem}[Bearing-Based Network Localization]\label{problem_bearingbasednetworkLocalization}
Consider a network $\G(p)$ in $\R^d$, the bearing-based network localization problem is to determine the locations of the follower nodes, $\{p_i\}_{i \in \V_f}$, given the inter-neighbor bearings, $\{g_{ij}\}_{(i,j)\in\E}$, and the locations of the anchor nodes, $\{p_i\}_{i\in\V_a}$.
Mathematically, the problem is to retrieve the true network location $p$ by solving the system of nonlinear equations,
\begin{align}\label{eq_networkLocalization_nonlinearConstraint}
\left\{
  \begin{array}{ll}
    \displaystyle{\frac{\hat{p}_j-\hat{p}_i}{\|\hat{p}_j-\hat{p}_i\|}=g_{ij}}, & \forall (i,j)\in\E, \\
    \hat{p}_i=p_i, & \forall i\in\V_a, \\
  \end{array}
\right.
\end{align}
where $\hat{p}_i$ is the estimated location of node $i$.
\end{problem}

The true network location is always a solution to the nonlinear equations in \eqref{eq_networkLocalization_nonlinearConstraint}, but the nonlinear equations may admit many other solutions that do not correspond to the true network location.
Thus we need to study when the true network location is the \emph{unique} solution to \eqref{eq_networkLocalization_nonlinearConstraint}, which motivates the following notion.

\begin{figure}[t]
  \centering
  \def\myscale{0.45}
  \def\mycolor{white}
\subfloat[]{
\begin{tikzpicture}[scale=\myscale]
            \def\length{3}
            \coordinate (x1) at (0,0);
            \coordinate (x4) at (\length,0);
            \coordinate (x3) at (\length,-\length);
            \coordinate (x2) at (0,-\length);
            \draw [] (x1)--(x2)--(x3)--(x4)--cycle;
            \def\radius{9pt}
            \draw [fill=black](x1) circle [radius=\radius];
            \draw (x1) node[] {\scriptsize{\color{white}$1$}};
            \draw [fill=black](x2) circle [radius=\radius];
            \draw (x2) node[] {\scriptsize{\color{white}$2$}};
            \draw [fill=white](x3) circle [radius=\radius];
            \draw (x3) node[] {\scriptsize{\color{black}$3$}};
            \draw [fill=white](x4) circle [radius=\radius];
            \draw (x4) node[] {\scriptsize{\color{black}$4$}};
\def\r{0.5}
\def\originX{0}
\def\originY{-\length/2}
\coordinate (x_left) at (\originX-\r,\originY);
\coordinate (x_right) at (\originX+\r,\originY);
\coordinate (y_up) at (\originX,\originY+\r);
\coordinate (y_down) at (\originX,\originY-\r);
\draw [blue](\originX,\originY) circle [radius=2pt];
\draw [->,thick,blue] (x_left)--(x_right)node[right=-2pt]{\tiny{$x$}};
\draw [->,thick,blue] (y_down)--(y_up) node[right=-1pt]{\tiny{$y$}};
\end{tikzpicture}}
\quad
\subfloat[]{
\begin{tikzpicture}[scale=\myscale]
            \def\length{3}
            \coordinate (x1) at (0,0);
            \coordinate (x4) at (\length/3*2,0);
            \coordinate (x3) at (\length/3*2,-\length);
            \coordinate (x2) at (0,-\length);
            \draw [] (x1)--(x2)--(x3)--(x4)--cycle;
            \def\radius{9pt}
            \draw [fill=black](x1) circle [radius=\radius];
            \draw (x1) node[] {\scriptsize{\color{white}$1$}};
            \draw [fill=black](x2) circle [radius=\radius];
            \draw (x2) node[] {\scriptsize{\color{white}$2$}};
            \draw [fill=white](x3) circle [radius=\radius];
            \draw (x3) node[] {\scriptsize{\color{black}$3$}};
            \draw [fill=white](x4) circle [radius=\radius];
            \draw (x4) node[] {\scriptsize{\color{black}$4$}};
\def\r{0.5}
\def\originX{0}
\def\originY{-\length/2}
\coordinate (x_left) at (\originX-\r,\originY);
\coordinate (x_right) at (\originX+\r,\originY);
\coordinate (y_up) at (\originX,\originY+\r);
\coordinate (y_down) at (\originX,\originY-\r);
\draw [blue](\originX,\originY) circle [radius=2pt];
\draw [->,thick,blue] (x_left)--(x_right)node[right=-2pt]{\tiny{$x$}};
\draw [->,thick,blue] (y_down)--(y_up) node[right=-1pt]{\tiny{$y$}};
\end{tikzpicture}}
\quad
\subfloat[]{
\begin{tikzpicture}[scale=\myscale]
            \def\length{3}
            \coordinate (x1) at (0,0);
            \coordinate (x4) at (-\length/3*2,0);
            \coordinate (x3) at (-\length/3*2,-\length);
            \coordinate (x2) at (0,-\length);
            \draw [] (x1)--(x2)--(x3)--(x4)--cycle;
            \def\radius{9pt}
            \draw [fill=black](x1) circle [radius=\radius];
            \draw (x1) node[] {\scriptsize{\color{white}$1$}};
            \draw [fill=black](x2) circle [radius=\radius];
            \draw (x2) node[] {\scriptsize{\color{white}$2$}};
            \draw [fill=white](x3) circle [radius=\radius];
            \draw (x3) node[] {\scriptsize{\color{black}$3$}};
            \draw [fill=white](x4) circle [radius=\radius];
            \draw (x4) node[] {\scriptsize{\color{black}$4$}};
\def\r{0.5}
\def\originX{0}
\def\originY{-\length/2}
\coordinate (x_left) at (\originX-\r,\originY);
\coordinate (x_right) at (\originX+\r,\originY);
\coordinate (y_up) at (\originX,\originY+\r);
\coordinate (y_down) at (\originX,\originY-\r);
\draw [blue](\originX,\originY) circle [radius=2pt];
\draw [->,thick,blue] (x_left)--(x_right)node[right=-2pt]{\tiny{$x$}};
\draw [->,thick,blue] (y_down)--(y_up) node[right=-1pt]{\tiny{$y$}};
\end{tikzpicture}}
  \caption{An illustration of the notion of localizability.
  Black dots represent the anchors and white dots for the followers.
  Suppose the true network is (a).
  The networks in (a) and (b) both satisfy the nonlinear equations in \eqref{eq_networkLocalization_nonlinearConstraint}.
  The networks in (a), (b), and (c) all satisfy the linear equations in \eqref{eq_networkLocalization_linearConstraint}.}
  \label{fig_example_illustrateLocalizability}
\end{figure}
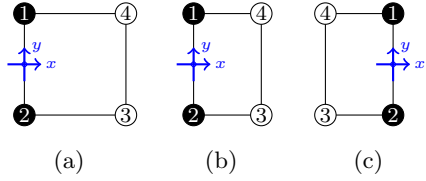

\begin{definition}[Bearing-Based Network Localizability]\label{definition_networkLocalizability}
A network $\G(p)$ is called \emph{bearing-based localizable} if the true network location $p$ is the unique solution to \eqref{eq_networkLocalization_nonlinearConstraint}.
\end{definition}

Localizability is a fundamental property of bearing-based networks.
A network must be localizable in order to be localized with either distributed or centralized protocols.
The notion of localizability is illustrated by an example in Figure~\ref{fig_example_illustrateLocalizability}.
In this example, the network in Figure~\ref{fig_example_illustrateLocalizability}(a) is the true network.
The network in Figure~\ref{fig_example_illustrateLocalizability}(b) has the same bearings and anchor locations as the true network.
As a result, both of the networks in Figure~\ref{fig_example_illustrateLocalizability}(a)-(b) are solutions to \eqref{eq_networkLocalization_nonlinearConstraint} and hence the networks are \emph{not} localizable by Definition~\ref{definition_networkLocalizability}.

For the sake of simplicity, we assume that the graph $\G$ is \emph{undirected}, which means $(i,j)\in\E\Leftrightarrow (j,i)\in\E$.
If the graph is directed, suppose $(i,j)\in\E$ but $(j,i)\notin\E$.
We can always add the edge $(j,i)$ into $\E$ to convert the directed graph to an undirected one.
The directed edges $(i,j)$ and $(j,i)$ imply two equations ${(\hat{p}_j-\hat{p}_i)}/{\|\hat{p}_j-\hat{p}_i\|}=g_{ij}$ and ${(\hat{p}_i-\hat{p}_j)}/{\|\hat{p}_i-\hat{p}_j\|}=g_{ji}$, respectively.
The two equations are equivalent because $g_{ji}=-g_{ij}$.
As a result, adding the edge $(j,i)$ does not affect the solutions to \eqref{eq_networkLocalization_nonlinearConstraint}.

\subsection{Reformulation as a Least-Squares Problem}
In order solve the nonlinear equations in \eqref{eq_networkLocalization_nonlinearConstraint}, we derive a companion system of linear equations.
In this direction, we first introduce a useful orthogonal projection operator.
For any nonzero vector $x\in\R^d$ ($d\ge2$), define the orthogonal projection operator $P: \R^d\rightarrow\R^{d\times d}$ as
\begin{align*}
    P(x) \triangleq I_d - \frac{x}{\|x\|}\frac{x^\T }{\|x\|}.
\end{align*}
For notational simplicity, denote $P_x \triangleq P(x)$.
The matrix $P_x$ geometrically projects any vector onto the orthogonal compliment of $x$.
It can be easily verified that $P_x^\T =P_x$, $P_x^2=P_x$, $\Null(P_x)=\myspan{x}$, and the eigenvalues of $P_x$ are $\{0,1^{(d-1)}\}$.

Consider now the projection matrix, $P_{g_{ij}}= I_d-g_{ij}g_{ij}^\T$, associated with the bearing $g_{ij}$.  By multiplying $P_{g_{ij}}$ on both sides of the first equation in \eqref{eq_networkLocalization_nonlinearConstraint}, the nonlinear algebraic problem \eqref{eq_networkLocalization_nonlinearConstraint} is converted to a system of linear equations,
\begin{align}\label{eq_networkLocalization_linearConstraint}
\left\{
  \begin{array}{ll}
    P_{g_{ij}}(\hat{p}_j-\hat{p}_i)=0, & \forall (i,j)\in\E, \\
    \hat{p}_i=p_i, & \forall i\in\V_a. \\
  \end{array}
\right.
\end{align}
System \eqref{eq_networkLocalization_linearConstraint} is {not} equivalent to system \eqref{eq_networkLocalization_nonlinearConstraint} in general.
But we have the exact relation between \eqref{eq_networkLocalization_nonlinearConstraint} and \eqref{eq_networkLocalization_linearConstraint} as described in the following lemma.

\newcommand{\X}{\mathcal{X}}
\begin{lemma}\label{lemma_nonlinearLinearConstraintsEquivalent}
Let $\X_1$ and $\X_2$ denote the set of all solutions satisfying \eqref{eq_networkLocalization_nonlinearConstraint} and \eqref{eq_networkLocalization_linearConstraint}, respectively.  Then
\begin{enumerate}[(a)]
\item $\{p\}\subseteq\X_1 \subseteq \X_2$;
\item $\{p\}=\X_1$ if and only if $\{p\}=\X_2$.
\end{enumerate}
\end{lemma}
\begin{pf}
(a)~Since the true network location $p$ is always a solution to \eqref{eq_networkLocalization_nonlinearConstraint} and  \eqref{eq_networkLocalization_linearConstraint}, we know $\X_1$ and $\X_2$ are nonempty and $\{p\}\subseteq\X_1$ and $\{p\} \subseteq \X_2$.
Since \eqref{eq_networkLocalization_linearConstraint} is obtained by multiplying  \eqref{eq_networkLocalization_nonlinearConstraint} by $P_{g_{ij}}$, we know any solution to \eqref{eq_networkLocalization_nonlinearConstraint} is also a solution to \eqref{eq_networkLocalization_linearConstraint}, showing $\X_1\subseteq \X_2$.

(b)~(\emph{Sufficiency}) Suppose $\{p\}=\X_2$.
It then follows from $\{p\}\subseteq\X_1 \subseteq \X_2$ that $\{p\}=\X_1$.
(\emph{Necessity})
Suppose $\{p\}=\X_1$.
We next prove $\{p\}=\X_2$ by contradiction.
Assume $p'\in\X_2$ and $p'\ne p$.
Let $\delta p\triangleq p'-p$ and define
\begin{align}\label{eq_definitionp''}
p''\triangleq p+k\delta p, \quad k\in\R.
\end{align}
We next show that $p''\in\X_1$ when $|k|$ is sufficiently small, leading to a contradiction.
Since $p, p'\in\X_2$, we know $p''\in\X_2$ for all $k\in\R$ by \eqref{eq_definitionp''}.
As a result, for any $k\in\R$ and $(i,j)\in\E$, we have $P_{{g}_{ij}}(p''_j-p''_i)=0$ which implies either $(p''_j-p''_i)/\|p''_j-p''_i\|=g_{ij}$ or $(p''_j-p''_i)/\|p''_j-p''_i\|=-g_{ij}$.
Since $p''_j-p''_i=(p_j-p_i)+k(\delta p_j-\delta p_i)$ according to \eqref{eq_definitionp''}, it is obvious that when $|k|$ is sufficiently small, the entries of $p''_j-p''_i$ have the same signs as those of $p_j-p_i$, and consequently $(p''_j-p''_i)/\|p''_j-p''_i\|=(p_j-p_i)/\|p_j-p_i\|=g_{ij}$.
Note that when any entry of $p_j-p_i$ is zero, the corresponding entry of $\delta p_j-\delta p_i$ is also zero because $\delta p_i-\delta p_j$ is parallel to $p_j-p_i$.
To conclude, $p''$ is another solution other than $p$ satisfying \eqref{eq_networkLocalization_nonlinearConstraint}, which is a contradiction.
\qed\end{pf}
\begin{remark}
The proof of Lemma~\ref{lemma_nonlinearLinearConstraintsEquivalent}(b) can be illustrated by Figure~\ref{fig_example_illustrateLocalizability}, where the networks (a), (b), and (c) correspond to $p$, $p''$, and $p'$ in the proof, respectively.
\end{remark}

Lemma~\ref{lemma_nonlinearLinearConstraintsEquivalent} indicates that the true network location $p$ is the unique solution to \eqref{eq_networkLocalization_nonlinearConstraint} if and only if $p$ is the unique solution to \eqref{eq_networkLocalization_linearConstraint}.
Thus we can study the localizability by analyzing the linear system \eqref{eq_networkLocalization_linearConstraint}.
The linear system of equations in \eqref{eq_networkLocalization_linearConstraint} can be rewritten as the following linear least-squares problem,
\begin{align}\label{eq_networkLocalizationLSOptimization}
  \underset{\hat{p}\in\R^{dn}}{\text{minimize}} & \qquad J(\hat{p})=\frac{1}{2}\sum_{i\in\V}\sum_{j\in\N_i}\|P_{g_{ij}}(\hat{p}_i-\hat{p}_j)\|^2, \\
    \text{subject to} &\qquad \hat{p}_i=p_i, \quad i\in\V_a.\nonumber
\end{align}
Since any minimizer with the objective function as zero is the solution to \eqref{eq_networkLocalization_linearConstraint}, we now successfully formulate the localizability problem as the above least-squares problem.
The rest of the paper is dedicated to studying two properties of the least-squares problem.
The first is to determine when the true location $p$ is the unique global minimizer of \eqref{eq_networkLocalizationLSOptimization} (i.e., the network is localizable), and the second is how to obtain $p$ in a distributed manner.

\section{The Bearing Laplacian Matrix}\label{section_propertiesOfBearingLaplacian}

In this section, we show that a new important matrix, termed \emph{bearing Laplacian}, emerges in the least-squares formulation.
The useful properties of the bearing Laplacian that will be used throughout the paper are explored.

Since the underlying graph $\G$ is undirected, the objective function in \eqref{eq_networkLocalizationLSOptimization} can be expressed in a quadratic form,
\begin{align*}
    J(\hat{p})
    =\hat{p}^\T \B(\G(p)) \hat{p},
\end{align*}
where $\B(\G(p))\in\R^{dn\times dn}$ and its $ij$th subblock matrix is
\begin{align*}
[\B(\G(p))]_{ij}=\left\{
  \begin{array}{ll}
      \mathbf{0}_{d\times d}, & i\ne j, (i,j)\notin\E, \\
      -P_{g_{ij}}, & i\ne j, (i,j)\in\E, \\ 
      \sum_{k\in\N_i}P_{g_{ik}}, & i=j, i\in\V. \\
  \end{array}
\right.
\end{align*}
For notational simplicity, we write $\B(\G(p))$ as $\B$ in the sequel.
The matrix $\B$ has a structure reminiscent of the weighted graph Laplacian matrix.
Since $\B$ indicates not only the topology of the network but also the inter-neighbor bearings, it is referred to as \emph{bearing Laplacian} in this paper.

The bearing Laplacian has an intimate connection to the \emph{bearing rigidity} properties of the network.
Preliminaries to the bearing rigidity theory, originally proposed in \cite{zhao2014TACBearing}, are given in Appendix~\ref{appendix_preliminaryBearingRigidity}.
Here we would like to highlight two important notions from this theory.
The first is the notion of \emph{infinitesimal bearing motions}.
Loosely speaking, infinitesimal bearing motions are motions of the nodes that preserve inter-neighbor bearings.
For example, for the network in Figure~\ref{fig_example_illustrateLocalizability}(a), the bearings can be preserved when the nodes 3 and 4 move in the horizontal direction to the right.
A network always has two kinds of \emph{trivial} infinitesimal bearing motions - they are the translational and scaling motions of the entire network.
A network is \emph{infinitesimally bearing rigid} if all its infinitesimal bearing motions are trivial.
One important property of an infinitesimally bearing rigid network is that its shape can be uniquely determined by the inter-neighbor bearings.

We next give the basic properties of the bearing Laplacian matrix.
We also show that the bearing Laplacian matrix is a powerful tool for characterizing the bearing rigidity of a network.

\begin{lemma}\label{lemma_NetworkLaplacian_nullspace}
    For a network $\G(p)$ with undirected graph $\G$, the bearing Laplacian $\B$ satisfies the following:
     \begin{enumerate}[(a)]
     \item $\B$ is symmetric positive semi-definite;
     \item $\rank(\B)\le dn-d-1$ and $\Null(\B)\supseteq\myspan{\one\otimes I_d, p}$;
     \item $\rank(\B)=dn-d-1$ and $\Null(\B)=\myspan{\one\otimes I_d, p}$ if and only if $\G(p)$ is infinitesimally bearing rigid.
     \end{enumerate}
\end{lemma}
\begin{pf}
Assign an arbitrary orientation to each undirected edge and label the edge vectors and bearings for the directed edges as $\{e_{k}\}_{k=1}^m$ and $\{g_{k}\}_{k=1}^m$, respectively.
Then the bearing Laplacian $\B$ can be expressed as $\B=\bar{H}^\T \mydiag(P_{g_k})\bar{H}$ where $\bar{H}=H\otimes I_d$ and $H$ is the incidence matrix of the graph.\footnote{The incidence matrix $H \in \R^{m\times n}$ of an oriented graph is the $\{0,\pm 1\}$-matrix with $[H]_{ki}=1$ if vertex $i$ is the head of edge $k$, $[H]_{ki}=-1$ if it is the tail, and $0$ otherwise.}
It further follows from $P_{g_k}=P_{g_k}^\T P_{g_k}$ that
    \begin{align*}
        \B=\underbrace{\bar{H}^\T \mydiag{(P_{g_k}^\T )}}_{\mathcal{R}^\T }\underbrace{\mydiag{(P_{g_k})}\bar{H}}_{\mathcal{R}}
        =\mathcal{R}^\T \mathcal{R}.
    \end{align*}
    Note $\mathcal{R}=\dia{\|e_k\|I_d}R_B$ where $R_B$ is the bearing rigidity matrix (see Lemma~\ref{lemma_bearingRigidityMatrixRank} in Appendix~\ref{appendix_preliminaryBearingRigidity}).
    As a result, the matrix $\mathcal{R}$, and hence $\B$, have exactly the same rank and null space as $R_B$.
    Then the results in (b) and (c) follows immediately from Lemma~\ref{lemma_bearingRigidityMatrixRank} and Theorem~\ref{theorem_conditionInfiParaRigid} as given in Appendix~\ref{appendix_preliminaryBearingRigidity}.
\qed\end{pf}

Since the nodes in the network are partitioned into anchors and followers, it will be useful to partition the corresponding bearing Laplacian as
\begin{align*}
    \B=\left[
         \begin{array}{cc}
           \B_{aa} & \B_{af} \\
           \B_{fa} & \B_{ff} \\
         \end{array}
       \right],
\end{align*}
where $\B_{aa}\in\R^{dn_a\times dn_a}$, $\B_{af}=\B_{fa}^\T \in\R^{dn_a\times dn_f}$, and $\B_{ff}\in\R^{dn_f\times dn_f}$.

\begin{lemma}\label{lemma_pa_pf_relation}
For any network $\G(p)$ with undirected graph $\G$, the subblock matrix $\B_{ff}$ is symmetric positive semi-definite and satisfies $\B_{ff}p_f+\B_{fa}p_a=0$.
\end{lemma}

\begin{pf}
For any nonzero $x\in\R^{dn_f}$, denote $\bar{x}=[0,x^\T ]^\T \in\R^{dn}$. Since $\B\ge0$, we have $x^\T \B_{ff}x=\bar{x}^\T \B\bar{x}\ge0$.
As a result $\B_{ff}$ is positive semi-definite.
Since $p\in\Null(\B)$ as suggested by Lemma~\ref{lemma_NetworkLaplacian_nullspace}, we have $\B p=0$ which further implies $\B_{fa}p_a+\B_{ff}p_f=0$.
\qed\end{pf}

\section{Analysis of Network Localizability}\label{section_localizabilityAnalysis}

In this section, we analyze the localizability of networks in arbitrary dimensions.
We first prove two necessary and sufficient conditions for network localizability from algebraic and rigidity perspectives, respectively.
We then present more necessary and/or sufficient conditions which can give more intuition on what localizable networks look like.

First of all, we derive the optimality condition for the least-squares problem \eqref{eq_networkLocalizationLSOptimization}.

\begin{lemma}\label{lemma_LSOptimalityCondition}
For the least-squares problem \eqref{eq_networkLocalizationLSOptimization}, any minimizer $\hat{p}_f^*$ is also a global minimizer and satisfies $$\B_{ff}\hat{p}_f^*+\B_{fa}p_a=0.$$
\end{lemma}
\begin{pf}
By substituting $\hat{p}_a=p_a$ into the objective function $J(\hat{p})=\hat{p}^\T \B \hat{p}$, the constrained optimization problem \eqref{eq_networkLocalizationLSOptimization} can be converted to the unconstrained problem
\begin{align}\label{eq_networkLocalizationLSOptimization_unconstrained}
    \min_{\hat{p}_f\in\R^{dn_f}}
    \tilde{J}(\hat{p}_f)=\hat{p}_f^\T \B_{ff}\hat{p}_f+2{p}_a^\T \B_{af}\hat{p}_f+p_a^\T \B_{aa}p_a.
\end{align}
Any minimizer must satisfy $\nabla_{\hat{p}_f}\tilde{J}(\hat{p}_f)=\B_{ff}\hat{p}_f+\B_{fa}p_a=0$.
Now suppose $\hat{p}_f^*$ is a minimizer and satisfies $\B_{ff}\hat{p}_f^*+\B_{fa}p_a=0$.
By comparing with $\B_{ff}p_f+\B_{fa}p_a=0$ as shown in Lemma~\ref{lemma_pa_pf_relation}, we know $\hat{p}_f^*=p_f+x$ where $x\in\Null(\B_{ff})$.
Let $\hat{p}^*=[p_a^\T,(\hat{p}_f^*)^\T]^\T$ and $\bar{x}=[0,x^\T]^\T\in\R^{dn}$. Since $\hat{p}_f^*={p}_f+x$ and $\B p=0$, we have $J(\hat{p}^*)=(\hat{p}^*)^\T\B\hat{p}^*=(p+\bar{x})^\T\B(p+\bar{x})=\bar{x}^\T\B\bar{x}=x^\T\B_{ff}x=0$.
As a result, the objective function equals zero at every minimizer.
\qed\end{pf}

The linear equations in \eqref{eq_networkLocalization_linearConstraint} hold if and only if the objective function in the least-squares problem \eqref{eq_networkLocalizationLSOptimization} is minimized to zero; this is a direct consequence of the first-order optimality conditions associated with \eqref{eq_networkLocalizationLSOptimization}.
Thus the equivalence between \eqref{eq_networkLocalization_linearConstraint} and \eqref{eq_networkLocalizationLSOptimization} is formally established.
We are now ready to present the necessary and sufficient condition for localizability.

\begin{theorem}[Algebraic Condition for Localizability]\label{theorem_Localizability_NS_alge}
A network $\G(p)$ is localizable if and only if the matrix $\B_{ff}$ is nonsingular.
When the network is localizable, the true locations of the followers can be calculated by $p_f=-\B_{ff}^{-1}\B_{fa}p_a$.
\end{theorem}
\begin{pf}
By Lemma~\ref{lemma_LSOptimalityCondition}, a network is localizable if and only if the true network location $p$ is the unique minimizer of the least-squares problem \eqref{eq_networkLocalizationLSOptimization}.
Since any minimizer must satisfy $\B_{ff}\hat{p}_f^*+\B_{fa}p_a=0$, it is obvious that the minimizer is unique if and only if $\B_{ff}$ is nonsingular.
When $\B_{ff}$ is nonsingular, we have $\hat{p}_f^*=-\B_{ff}^{-1}\B_{fa}p_a$, whose value equals the true location $p_f$ according to Lemma~\ref{lemma_pa_pf_relation}.
\qed\end{pf}

Theorem~\ref{theorem_Localizability_NS_alge} establishes the equivalence between the localizability and the nonsingularity of $\B_{ff}$.
A question that immediately follows Theorem~\ref{theorem_Localizability_NS_alge} is what kind of networks have nonsingular $\B_{ff}$.
We next propose a necessary and sufficient condition from the bearing rigidity point of view.
This rigidity condition is mathematically equivalent to the algebraic condition, but it gives more intuition on what localizable networks look like.

\begin{theorem}[Rigidity Condition for Localizability]\label{theorem_Localizability_NS_rigidity}
A network $\G(p)$ is localizable if and only if every infinitesimal bearing motion involves at least one anchor; that is, for any nonzero infinitesimal bearing motion
\begin{align*}
    \delta p=\left[
               \begin{array}{c}
                 \delta p_a \\
                  \delta p_f \\
               \end{array}
             \right]\in\Null(\B),
\end{align*}
the vector $\delta p_a$ corresponding to the anchors must be nonzero.
\end{theorem}
\begin{pf}
We only need to show that $\B_{ff}$ is singular if and only if there exists nonzero $\delta p\in\Null(\B)$ with $\delta p_a=0$.
(\emph{Necessity}) Suppose $\B_{ff}$ is singular.
Then there exists nonzero $x\in\R^{dn_f}$ such that $\B_{ff}x=0$.
Let $\delta p=[0,x^\T]^\T\in\R^{dn}$. Then $\delta p^\T\B \delta p=x^\T\B_{ff}x=0$.
Hence $\delta p\in \Null(\B)$ and $\delta p_a=0$.
(\emph{Sufficiency}) Suppose there exists $\delta p\in\Null(\B)$ satisfying $\delta p_a=0$ and $\delta p_f\ne0$.
Then $\delta p_f^\T\B_{ff}\delta p_f=\delta p^\T \B \delta p=0$, which implies that $\B_{ff}$ is singular.
\qed\end{pf}

\begin{figure*}
  \centering
  \def\myscale{0.35}
  \subfloat[]{
\begin{tikzpicture}[scale=\myscale]
\def\arrowL{1}
\def\offset{6pt}
\def\length{5}
\coordinate (x1) at (0,0);
\coordinate (x1Above) at (0,\offset);
\coordinate (x1Below) at (0,-\offset);
\coordinate (x2) at (\length/2,\offset);
\coordinate (x3) at (\length,0);
\coordinate (x3Above) at (\length,\offset);
\coordinate (x3Below) at (\length,-\offset);
\draw [semithick] (x1Above)--(x2)--(x3Above);
\draw [semithick] (x1Below)--(x3Below);
\draw [thick,->,draw=red] (x2)--(\length/2+\arrowL,\offset);
\def\radius{6pt}
\draw [fill=black](x1) circle [radius=\radius];
\draw [fill=white](x2) circle [radius=\radius];
\draw [fill=black](x3) circle [radius=\radius];
\end{tikzpicture}}
\quad
\subfloat[]{
\begin{tikzpicture}[scale=\myscale]
\def\arrowL{0.8}
\def\length{4}
            \coordinate (x1) at (0,0);
            \coordinate (x2) at (4,0);
            \coordinate (x3) at (4,-4);
            \coordinate (x4) at (0,-4);
            \draw [semithick] (x1)--(x2)--(x3)--(x4)--cycle;
            \draw [thick,->,draw=red] (x1)--(0,\arrowL);
            \draw [thick,->,draw=red] (x2)--(\length,\arrowL);
            \def\radius{6pt}
            \draw [fill=white](x1) circle [radius=\radius];
            \draw [fill=white](x2) circle [radius=\radius];
            \draw [fill=black](x3) circle [radius=\radius];
            \draw [fill=black](x4) circle [radius=\radius];
\end{tikzpicture}}
\quad
\subfloat[]{
\begin{tikzpicture}[scale=\myscale]
\def\length{5}
\def\arrowL{0.8}
            \def\scaleBig{3}; 
            \coordinate (x1) at (-0.866*\scaleBig,-0.5*\scaleBig);
            \coordinate (x2) at (0.866*\scaleBig,-0.5*\scaleBig);
            \coordinate (x3) at (0,\scaleBig);
            \def\scaleSmall{1.5}; 
            \coordinate (x4) at (-0.866*\scaleSmall,-0.5*\scaleSmall);
            \coordinate (x5) at (0.866*\scaleSmall,-0.5*\scaleSmall);
            \coordinate (x6) at (0,\scaleSmall);
            \draw [semithick] (x1)--(x4);
            \draw [semithick] (x2)--(x5);
            \draw [semithick] (x3)--(x6);
            \draw [semithick] (x1)--(x2)--(x3)--cycle;
            \draw [semithick] (x4)--(x5)--(x6)--cycle;
            \draw [thick,->,draw=red] (x4)--(-0.866*\scaleSmall-0.866*\arrowL,-0.5*\scaleSmall-0.5*\arrowL);
            \draw [thick,->,draw=red] (x5)--(0.866*\scaleSmall+0.866*\arrowL,-0.5*\scaleSmall-0.5*\arrowL);
            \draw [thick,->,draw=red] (x6)--(0,\scaleSmall+\arrowL);
            \def\radius{6pt}
            \draw [fill=black](x1) circle [radius=\radius];
            \draw [fill=black](x2) circle [radius=\radius];
            \draw [fill=black](x3) circle [radius=\radius];
            \draw [fill=white](x4) circle [radius=\radius];
            \draw [fill=white](x5) circle [radius=\radius];
            \draw [fill=white](x6) circle [radius=\radius];
\end{tikzpicture}}
\quad
\subfloat[]{
\begin{tikzpicture}[scale=\myscale]
\def\arrowL{1}
\def\length{5}
\def\high{4}
\coordinate (x1) at (0,0);
\coordinate (x2) at (\length,0);
\coordinate (x3) at (\length,\high);
\coordinate (x4) at (0,\high);
\coordinate (x5) at (\length/3,\high/2);
\coordinate (x6) at (\length/3*2,\high/2);
\draw [semithick] (x1)--(x2)--(x3)--(x4)--cycle;
\draw [semithick] (x1)--(x5)--(x4);
\draw [semithick] (x2)--(x6)--(x3);
\draw [semithick] (x5)--(x6);
\draw [thick,->,draw=red] (x3)--(\length,\high+\arrowL);
\draw [thick,->,draw=red] (x4)--(0,\high+\arrowL);
\draw [thick,->,draw=red] (x5)--(\length/3+0.64*0.707*\arrowL,\high/2+0.7682*0.707*\arrowL);
\draw [thick,->,draw=red] (x6)--(\length/3*2-0.64*0.707*\arrowL,\high/2+0.7682*0.707*\arrowL);
\def\radius{6pt}
\draw [fill=black](x1) circle [radius=\radius];
\draw [fill=black](x2) circle [radius=\radius];
\draw [fill=white](x3) circle [radius=\radius];
\draw [fill=white](x4) circle [radius=\radius];
\draw [fill=white](x5) circle [radius=\radius];
\draw [fill=white](x6) circle [radius=\radius];
\end{tikzpicture}}
\quad
\subfloat[]{
\begin{tikzpicture}[scale=\myscale]
\def\arrowL{0.8}
            \def\length{2}
            \def\height{3.2}
            \coordinate (x1) at (0,-1/3*\height);
            \coordinate (x2) at (\length,0);
            \coordinate (x3) at (-\length,0);
            \coordinate (x4) at (0,-4/3*\height);
            \coordinate (x5) at (\length,-\height);
            \coordinate (x6) at (-\length,-\height);
            \draw [semithick] (x1)--(x2)--(x3)--cycle;
            \draw [semithick] (x6)--(x4)--(x5);
            \draw [semithick,densely dotted] (x5)--(x6);
            \draw [semithick] (x1)--(x4);
            \draw [semithick] (x2)--(x5);
            \draw [semithick] (x3)--(x6);
            \draw [thick,->,draw=red] (x1)--(0,-\height/3+\arrowL);
            \draw [thick,->,draw=red] (x2)--(\length,\arrowL);
            \draw [thick,->,draw=red] (x3)--(-\length,\arrowL);
             draw circles 
            \def\radius{6pt}
            \draw [fill=white](x1) circle [radius=\radius];
            \draw [fill=white](x2) circle [radius=\radius];
            \draw [fill=white](x3) circle [radius=\radius];
            \draw [fill=white](x4) circle [radius=\radius];
            \draw [fill=black](x5) circle [radius=\radius];
            \draw [fill=black](x6) circle [radius=\radius];
\end{tikzpicture}}
\quad
\subfloat[]{
\begin{tikzpicture}[scale=\myscale]
\def\arrowL{0.8}
            \def\length{3.5}
            \coordinate (x1) at (0,0);
            \coordinate (x2) at (\length,0);
            \coordinate (x3) at (\length,-\length);
            \coordinate (x4) at (0,-\length);
            \def\Xoffset{1}
            \def\Yoffset{1}
            \coordinate (x5) at (0+\Xoffset,0+\Yoffset);
            \coordinate (x6) at (\length+\Xoffset,0+\Yoffset);
            \coordinate (x7) at (\length+\Xoffset,-\length+\Yoffset);
            \coordinate (x8) at (0+\Xoffset,-\length+\Yoffset);
            \draw [semithick] (x1)--(x2)--(x3)--(x4)--cycle;
            \draw [semithick] (x1)--(x5)--(x6)--(x2);
            \draw [semithick] (x6)--(x7)--(x3);
            \draw [semithick,densely dotted] (x5)--(x8);
            \draw [semithick,densely dotted] (x7)--(x8);
            \draw [semithick,densely dotted] (x4)--(x8);
            \draw [thick,->,draw=red] (x1)--(0,0+\arrowL);
            \draw [thick,->,draw=red] (x2)--(\length,0+\arrowL);
            \draw [thick,->,draw=red] (x5)--(\Xoffset,\Yoffset+\arrowL);
            \draw [thick,->,draw=red] (x6)--(\length+\Xoffset,\Yoffset+\arrowL);                                    
            \def\radius{6pt}
            \draw [fill=white](x1) circle [radius=\radius];
            \draw [fill=white](x2) circle [radius=\radius];
            \draw [fill=black](x3) circle [radius=\radius];
            \draw [fill=black](x4) circle [radius=\radius];
            \draw [fill=white](x5) circle [radius=\radius];
            \draw [fill=white](x6) circle [radius=\radius];
            \draw [fill=white](x7) circle [radius=\radius];
            \draw [fill=white,densely dotted](x8) circle [radius=\radius];
\end{tikzpicture}
} 
  \caption{Examples of \emph{non-localizable} networks.
  The solid and hollow dots represent the anchors and followers, respectively.
  The networks are not localizable because they have infinitesimal bearing motions that only correspond to the followers (see, for example, the red arrows).
  The networks in (e) and (f) are three-dimensional, and the rest are two-dimensional.}
  \label{fig_Example_nonlocalizable}
\end{figure*}
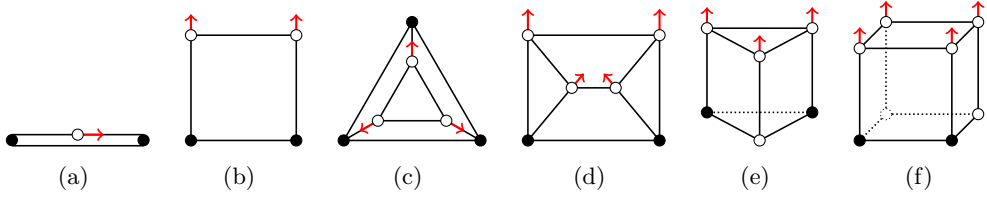
\begin{figure*}
  \centering
  \def\myscale{0.35}
\subfloat[]{
\begin{tikzpicture}[scale=\myscale]
            \coordinate (x1) at (0,0);
            \coordinate (x2) at (4,0);
            \coordinate (x3) at (2,3);
            \draw [semithick] (x1)--(x2)--(x3)--cycle;
            \def\radius{6pt}
            \draw [fill=black](x1) circle [radius=\radius];
            \draw [fill=black](x2) circle [radius=\radius];
            \draw [fill=white](x3) circle [radius=\radius];
\end{tikzpicture}}
\,
\subfloat[]{
\begin{tikzpicture}[scale=\myscale]
            \coordinate (x1) at (0,0);
            \coordinate (x2) at (4,0);
            \coordinate (x3) at (4,-4);
            \coordinate (x4) at (0,-4);
            \draw [semithick] (x1)--(x2)--(x3)--(x4)--cycle;
            \draw [semithick] (x2)--(x4);
            \def\radius{6pt}
            \draw [fill=white](x1) circle [radius=\radius];
            \draw [fill=white](x2) circle [radius=\radius];
            \draw [fill=black](x3) circle [radius=\radius];
            \draw [fill=black](x4) circle [radius=\radius];
\end{tikzpicture}}
\,
\subfloat[]{
\begin{tikzpicture}[scale=\myscale]
            \def\scaleBig{3}; 
            \coordinate (x1) at (-0.866*\scaleBig,-0.5*\scaleBig);
            \coordinate (x2) at (0.866*\scaleBig,-0.5*\scaleBig);
            \coordinate (x3) at (0,\scaleBig);
            \def\scaleSmall{1.5}; 
            \coordinate (x4) at (-0.866*\scaleSmall,-0.5*\scaleSmall);
            \coordinate (x5) at (0.866*\scaleSmall,-0.5*\scaleSmall);
            \coordinate (x6) at (0,\scaleSmall);
            \draw [semithick] (x1)--(x4);
            \draw [semithick] (x2)--(x5);
            \draw [semithick] (x3)--(x6);
            \draw [semithick] (x1)--(x2)--(x3)--cycle;
            \draw [semithick] (x4)--(x5)--(x6)--cycle;
            \def\radius{6pt}
            \draw [fill=black](x1) circle [radius=\radius];
            \draw [fill=black](x2) circle [radius=\radius];
            \draw [fill=white](x3) circle [radius=\radius];
            \draw [fill=white](x4) circle [radius=\radius];
            \draw [fill=white](x5) circle [radius=\radius];
            \draw [fill=black](x6) circle [radius=\radius];
\end{tikzpicture}}
\,
\subfloat[]{
\begin{tikzpicture}[scale=\myscale]
\def\length{5}
\def\high{4}
\coordinate (x1) at (0,0);
\coordinate (x2) at (\length,0);
\coordinate (x3) at (\length,\high);
\coordinate (x4) at (0,\high);
\coordinate (x5) at (\length/3,\high/2);
\coordinate (x6) at (\length/3*2,\high/2);
\draw [semithick] (x1)--(x2)--(x3)--(x4)--cycle;
\draw [semithick] (x1)--(x5)--(x4);
\draw [semithick] (x2)--(x6)--(x3);
\draw [semithick] (x5)--(x6);
\def\radius{6pt}
\draw [fill=black](x1) circle [radius=\radius];
\draw [fill=white](x2) circle [radius=\radius];
\draw [fill=white](x3) circle [radius=\radius];
\draw [fill=white](x4) circle [radius=\radius];
\draw [fill=white](x5) circle [radius=\radius];
\draw [fill=black](x6) circle [radius=\radius];
\end{tikzpicture}}
\,
\subfloat[]{
\begin{tikzpicture}[scale=\myscale]
\def\arrowL{0.8}
            \def\length{2}
            \def\height{3.2}
            \coordinate (x1) at (0,-1/3*\height);
            \coordinate (x2) at (\length,0);
            \coordinate (x3) at (-\length,0);
            \coordinate (x4) at (0,-4/3*\height);
            \coordinate (x5) at (\length,-\height);
            \coordinate (x6) at (-\length,-\height);
            \draw [semithick] (x1)--(x2)--(x3)--cycle;
            \draw [semithick] (x6)--(x4)--(x5);
            \draw [semithick,densely dotted] (x5)--(x6);
            \draw [semithick] (x1)--(x4);
            \draw [semithick] (x2)--(x5);
            \draw [semithick] (x3)--(x6);
             draw circles 
            \def\radius{6pt}
            \draw [fill=white](x1) circle [radius=\radius];
            \draw [fill=black](x2) circle [radius=\radius];
            \draw [fill=white](x3) circle [radius=\radius];
            \draw [fill=white](x4) circle [radius=\radius];
            \draw [fill=white](x5) circle [radius=\radius];
            \draw [fill=black](x6) circle [radius=\radius];
\end{tikzpicture}
}
\,
\subfloat[]{
\begin{tikzpicture}[scale=\myscale]
            \def\length{3.5}
            \coordinate (x1) at (0,0);
            \coordinate (x2) at (\length,0);
            \coordinate (x3) at (\length,-\length);
            \coordinate (x4) at (0,-\length);
            \def\Xoffset{1}
            \def\Yoffset{1}
            \coordinate (x5) at (0+\Xoffset,0+\Yoffset);
            \coordinate (x6) at (\length+\Xoffset,0+\Yoffset);
            \coordinate (x7) at (\length+\Xoffset,-\length+\Yoffset);
            \coordinate (x8) at (0+\Xoffset,-\length+\Yoffset);
            \draw [semithick] (x1)--(x2)--(x3)--(x4)--cycle;
            \draw [semithick] (x1)--(x5)--(x6)--(x2);
            \draw [semithick] (x6)--(x7)--(x3);
            \draw [semithick,densely dotted] (x5)--(x8);
            \draw [semithick,densely dotted] (x7)--(x8);
            \draw [semithick,densely dotted] (x4)--(x8);
            \def\radius{6pt}
            \draw [fill=white](x1) circle [radius=\radius];
            \draw [fill=white](x2) circle [radius=\radius];
            \draw [fill=white](x3) circle [radius=\radius];
            \draw [fill=black](x4) circle [radius=\radius];
            \draw [fill=white](x5) circle [radius=\radius];
            \draw [fill=black](x6) circle [radius=\radius];
            \draw [fill=white](x7) circle [radius=\radius];
            \draw [fill=white,densely dotted](x8) circle [radius=\radius];
\end{tikzpicture}
}
\,
\subfloat[]{
\begin{tikzpicture}[scale=\myscale]
            \coordinate (x1) at (1,0);
            \coordinate (x2) at (4,0);
            \coordinate (x3) at (4,4);
            \coordinate (x4) at (8,0);
            \coordinate (x6) at (8,4);
            \draw [semithick] (x1)--(x3)--(x2)--cycle;
            \draw [semithick] (x2)--(x4)--(x6)--(x3);
            \def\radius{6pt}
            \draw [fill=black](x1) circle [radius=\radius];
            \draw [fill=black](x2) circle [radius=\radius];
            \draw [fill=white](x3) circle [radius=\radius];
            \draw [fill=black](x4) circle [radius=\radius];
            \draw [fill=white](x6) circle [radius=\radius];
\end{tikzpicture}}
\,\subfloat[]{
\begin{tikzpicture}[scale=\myscale]
            \def\length{1.2}
            \coordinate (x1) at (-\length,-2/3*\length);
            \coordinate (x2) at (\length,-2/3*\length);
            \coordinate (x3) at (2*\length,0);
            \coordinate (x4) at (\length-0.05,2/3*\length-0.1);
            \coordinate (x5) at (-\length+0.05,2/3*\length-0.1);
            \coordinate (x6) at (-2*\length,0);
            \coordinate (x7) at (0,3*\length);
            \draw [semithick] (x6)--(x1)--(x2)--(x3);
            \draw [semithick,densely dotted] (x3)--(x4)--(x5)--(x6);
            \draw [semithick] (x1)--(x7);
            \draw [semithick] (x2)--(x7);
            \draw [semithick] (x3)--(x7);
            \draw [semithick,densely dotted] (x4)--(x7);
            \draw [semithick,densely dotted] (x5)--(x7);
            \draw [semithick] (x6)--(x7);
            \def\radius{6pt}
            \draw [fill=black](x1) circle [radius=\radius];
            \draw [fill=black](x2) circle [radius=\radius];
            \draw [fill=white](x3) circle [radius=\radius];
            \draw [fill=white,densely dotted](x4) circle [radius=\radius];
            \draw [fill=white,densely dotted](x5) circle [radius=\radius];
            \draw [fill=white](x6) circle [radius=\radius];
            \draw [fill=white](x7) circle [radius=\radius];
\end{tikzpicture}
} 
  \caption{Examples of \emph{localizable} networks.
  The solid and hollow dots represent the anchors and followers, respectively.
  The networks in (e), (f), and (h) are three-dimensional, and the rest are two-dimensional.}
  \label{fig_Example_localizable}
\end{figure*}
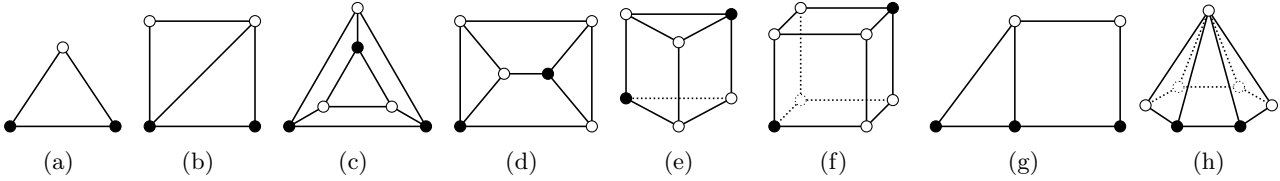

The intuition behind Theorem~\ref{theorem_Localizability_NS_rigidity} is as follows.
Any infinitesimal bearing motion (i.e., bearing-preserved motion) would imply multiple false networks that have exactly the same bearings as the true network.
Only if the infinitesimal bearing motion involves at least one anchor, the false networks can be ruled out as solutions to \eqref{eq_networkLocalization_nonlinearConstraint} since they do not satisfy the anchor constraints; otherwise, the false networks cannot be distinguished from the true network.

Examples are given in Figure~\ref{fig_Example_nonlocalizable} and Figure~\ref{fig_Example_localizable} to illustrate Theorem~\ref{theorem_Localizability_NS_rigidity}.
Figure~\ref{fig_Example_nonlocalizable} shows examples of \emph{non-localizable} networks.
These networks are not localizable because each of them has infinitesimal bearing motions that only involve the followers (see those marked by red arrows).
Figure~\ref{fig_Example_localizable} shows examples of \emph{localizable} networks.
The networks in Figure~\ref{fig_Example_localizable}(a)-(f) are obtained by modifying the networks in Figure~\ref{fig_Example_nonlocalizable}, which suggests that a non-localizable network can be made localizable by adding extra edges or selecting different anchors.
It is worth noting that the networks in Figure~\ref{fig_Example_localizable}(c)-(g) are not infinitesimally bearing rigid yet they are localizable.
As a result, infinitesimal bearing rigidity is not necessary to guarantee localizability.

Up to this point, we have presented two necessary and sufficient localizability conditions.
One is the algebraic condition in Theorem~\ref{theorem_Localizability_NS_alge} and the other is the rigidity condition in Theorem~\ref{theorem_Localizability_NS_rigidity}.
We next utilize the two conditions to examine some specific problems more closely.
The first is how many anchors are required to ensure the localizability of a network.

\begin{corollary}\label{corollary_anchorNumber}
    If a network $\G(p)$ is localizable, then
    \begin{align*}
        n_a\ge \frac{\dim\left(\Null(\B)\right)}{d}>1.
    \end{align*}
\end{corollary}
\begin{pf}
Let $k=\dim\left(\Null(\B)\right)$ and $N\in\R^{dn\times k}$ be a basis matrix of $\Null(\B)$ which means $\Range(N)=\Null(\B)$.
Then any nonzero $\delta p\in\Null(\B)$ can be expressed as $\delta p=N x$, where $x\in\R^k, x\ne 0$.
Partition $N$ and express $Nx$ as
$
\delta p=N x
=\left[
\begin{array}{c}
N_ax \\
N_fx \\
\end{array}
\right]
$,
where $N_a\in\R^{dn_a\times k}$.
According to Theorem~\ref{theorem_Localizability_NS_rigidity}, the network is localizable if and only if $N_ax\ne0, \forall x\in\R^k, x\ne 0$.
As a result, the matrix $N_a$ must have full column rank, which requires $N_a$ to be a \emph{tall} matrix with $dn_a\ge k=\dim(\Null(\B))$.
Since $\dim(\Null(\B))\ge d+1$ according to Lemma~\ref{lemma_NetworkLaplacian_nullspace}, we have $n_a\ge \dim(\Null(\B))/d\ge(d+1)/d>1$.
\qed\end{pf}

A simple but important fact suggested by Corollary~\ref{corollary_anchorNumber} is that any localizable network must have at least \emph{two} anchors.
Similar conclusions have already been obtained in the existing studies for networks in the two-dimensional space \cite{ErenTurkishJournal2007,Piovan2013Automatica,Shames2013TAC,ZhuGuangwei2014Automatica}.
But Corollary~\ref{corollary_anchorNumber} also suggests another important fact, which has not been observed in the literature, that more anchors are required to ensure the localizability when $\dim(\Null(\B))$ increases.
The quantity $\dim(\Null(\B))$ can be viewed as a measure of the ``degree of bearing rigidity'' because $\dim(\Null(\B))$ reaches the smallest value $d+1$ when the network is infinitesimally bearing rigid as shown in Lemma~\ref{lemma_NetworkLaplacian_nullspace}.
As a result, the intuition behind the second fact is that more anchors are required to ensure the localizability when the network is ``less'' bearing rigid (i.e., $\dim(\Null(\B))$ is large).

We next present another three localizability conditions, two of which are sufficient and the other is both necessary and sufficient.
These conditions are important because they indicate the explicit connection between the localizability and infinitesimal bearing rigidity.
Before presenting the conditions, we need to first define the notion of augmented network.

\begin{definition}[Augmented Network]
Given a network $\G(p)$ with $\G=(\V, \E)$, denote by $\bar{\G}(p)$ an \emph{augmented network} with $\bar{\G}=(\V, \bar{\E})$ where $\bar{\E}=\E\cup\{(i,j): i,j\in\V_a\}$.
\end{definition}

The augmented network $\bar{\G}(p)$ is obtained from $\G(p)$ by connecting every pair of anchors. 
If the anchors are already connected in $\G(p)$, then $\bar{\G}(p)$ is the same as $\G(p)$.
It should be noted that adding or deleting the edge between any pair of anchors only changes $\B_{aa}$ but not $\B_{ff}$.
As a result, $\G(p)$ and $\bar{\G}(p)$ have exactly the same $\B_{ff}$ and hence they are localizable or nonlocalizable simultaneously.
The next two sufficient conditions connect the notions of localizability and infinitesimal bearing rigidity.

\begin{corollary}\label{corollary_networkLocalizability_GBarIBR}
    When $n_a\ge2$, if $\bar{\G}(p)$ is infinitesimally bearing rigid, then $\G(p)$ is localizable.
\end{corollary}
\begin{pf}
We will first use Theorem~\ref{theorem_Localizability_NS_rigidity} to prove the localizability of $\bar{\G}(p)$.
Then the localizability of $\G(p)$ immediately follows because $\G(p)$ and $\bar{\G}(p)$ have the same localizability.
Let $\bar{\B}$ be the bearing Laplacian for $\bar{\G}(p)$.
Since $\bar{\G}(p)$ is infinitesimally bearing rigid, we have $\Null(\bar{\B})=\myspan{\one\otimes I_d, p}$ by Lemma~\ref{lemma_NetworkLaplacian_nullspace}.
As a result, any infinitesimal bearing motion $\delta p\in\Null(\bar{\B})$ can be expressed as a linear combination of $\one\otimes I_d$ and $p$.
Since no two anchors collocate, there does not exist a linear combination of $\one\otimes I_d$ and $p$ leading to $\delta p_a=0$ if $n_a\ge2$.
Then $\bar{\G}(p)$ is localizable according to Theorem~\ref{theorem_Localizability_NS_rigidity}.
\qed\end{pf}

\begin{corollary}\label{corollary_networkLocalizability_GIBR}
    When $n_a\ge2$, if $\G(p)$ is infinitesimally bearing rigid, then $\G(p)$ is localizable.
\end{corollary}
\begin{pf}
Similar to Corollary~\ref{corollary_networkLocalizability_GBarIBR}.
\qed\end{pf}

The intuition behind Corollary~\ref{corollary_networkLocalizability_GIBR} is as follows.
If a network is infinitesimally bearing rigid, then it can be uniquely determined up to a translation and a scaling factor by the bearings.
Since the translational and scaling ambiguity can be further eliminated by the anchor constraints, the entire network can be fully determined and hence localizable.
It is notable that Corollary~\ref{corollary_networkLocalizability_GIBR} is more restrictive than Corollary~\ref{corollary_networkLocalizability_GBarIBR} because it requires $\G(p)$ to be infinitesimally bearing rigid whereas Corollary~\ref{corollary_networkLocalizability_GBarIBR} merely requires $\bar{\G}(p)$ to be.
To illustrate, for each of the networks as shown in Figure~\ref{fig_Example_localizable}(c)-(f), the augmented network $\bar{\G}(p)$ is infinitesimally bearing rigid but $\G(p)$ is not.
Then, these networks can be concluded as localizable by Corollary~\ref{corollary_networkLocalizability_GBarIBR}.
Finally, Corollary~\ref{corollary_networkLocalizability_GBarIBR} can be viewed as a generalization of the result \cite[Cor~10]{ZhuGuangwei2014Automatica} which is applicable only to two-dimensional cases.

As suggested by Corollary~\ref{corollary_networkLocalizability_GBarIBR}, the condition of the infinitesimal bearing rigidity of $\bar{\G}(p)$ is \emph{sufficient} to ensure the localizability of $\G(p)$.
An important yet unexplored problem is whether or not the condition is also \emph{necessary}.
In the case of $n_a\ge3$, the condition is sufficient but \emph{not} necessary.
For example, for the network in Figure~\ref{fig_Example_localizable}(g), $\G(p)$ is localizable but $\bar{\G}(p)$ is not infinitesimally bearing rigid since the three anchors are collinear.
However, in the case of $n_a=2$, the condition is \emph{both necessary and sufficient} as shown below.

\begin{theorem}\label{theorem_na=2_IBRisNS}
When $n_a=2$, a network $\G(p)$ is localizable if and only if the augmented network $\bar{\G}(p)$ is infinitesimal bearing rigid.
\end{theorem}
\begin{pf}
The sufficiency has already been proved in Corollary~\ref{corollary_networkLocalizability_GBarIBR}.
We next prove the necessity by contradiction.
Assume $\G(p)$ is localizable but $\bar{\G}(p)$ is \emph{not} infinitesimal bearing rigid.
Then $\bar{\G}(p)$ has a nontrivial infinitesimal bearing motion $\delta p$ which is not in $\myspan{\one\otimes I_d, p}$.
Write $\delta p=[\delta p_1^\T, \delta p_2^\T, (*)]^\T$, where $\delta p_1, \delta p_2\in\R^d$ corresponds to the two anchors.
Because the infinitesimal motion $\delta p$ preserves all the bearings including the bearing between $p_1$ and $p_2$, we know that the vector $\delta p_1-\delta p_2$ is parallel to $p_1-p_2$.
As a result, there exists a nonzero scalar $k$ such that $\delta p_1-\delta p_2=k(p_1-p_2)$. Construct
\begin{align*}
\delta p'
&\triangleq \delta p+\one_n\otimes (kp_2-\delta p_2) -kp \\
&= \left[
  \begin{array}{c}
    \delta p_1 \\
    \delta p_2 \\
    (*)\\
  \end{array}
\right]
+\left[\begin{array}{c}
    kp_2-\delta p_2 \\
    kp_2-\delta p_2 \\
    (*)\\
  \end{array}
\right]
-\left[\begin{array}{c}
    k p_1 \\
    k p_2 \\
    (*)\\
  \end{array}
\right]
=
\left[
  \begin{array}{c}
    0 \\
    0 \\
    (*)\\
  \end{array}
\right].
\end{align*}
Since the first two entries of $\delta p'$ are zero, we know $\delta p'$ is an infinitesimal motion that only involves the followers.
Thus, the network is not localizable by Theorem~\ref{theorem_Localizability_NS_rigidity}, which is a contradiction.
\qed\end{pf}

\section{Distributed Network Localization Protocols}\label{section_distributedProtocol}

In this section, we propose and analyze a linear distributed protocol for bearing-based network localization in arbitrary dimensions.

The global minimizer of the unconstrained optimization problem \eqref{eq_networkLocalizationLSOptimization_unconstrained} can be obtained by the gradient decent protocol
\begin{align}\label{eq_bearingNetworkLocalizeProtocal}
    \dot{\hat{p}}_f(t)
    =-\nabla_{\hat{p}_f}\tilde{J}(\hat{p}_f)
    =-\B_{ff}\hat{p}_f(t)-\B_{fa}p_a,
\end{align}
whose elementwise expression is
\begin{align}\label{eq_bearingNetworkLocalizeProtocal_element}
    \dot{\hat{p}}_i(t)&=-\sum_{j\in\N_i} P_{g_{ij}}(\hat{p}_i(t)-\hat{p}_j(t)), \quad i\in\V_f.
\end{align}%
where $P_{g_{ij}}=I_d-g_{ij}g_{ij}^\T$.
Note the neighbor of the follower $i$ can be either a follower or an anchor.

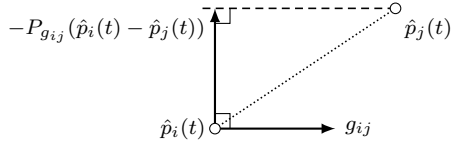
\begin{figure}
  \centering
  \def\myscale{0.4}
  \begin{tikzpicture}[scale=\myscale]
            \coordinate (xi) at (0,0);
            \coordinate (xj) at (6,4);
            \coordinate (proj) at (0,4);
            \coordinate (projLeft) at (-0.5,4);
            \def\unit{4}
            \coordinate (gij)                 at (\unit,0);
            \draw [densely dotted, semithick] (xi)--(xj);
            \draw [densely dashed, semithick] (xj)--(projLeft);
            \draw[->, >=latex, thick] (xi) -- (gij) node [right] {\scriptsize{$g_{ij}$}};
            \draw[->, >=latex, thick] (xi) -- (proj) node [below left] {\scriptsize{$-P_{{g}_{ij}}(\hat{p}_i(t)-\hat{p}_j(t))$}};
            \def\radius{5pt}
            \draw [fill=white](xi) circle [radius=\radius];
            \draw [fill=white](xj) circle [radius=\radius];
            \draw (xi) node[left] {\scriptsize{$\hat{p}_i(t)$}};
            \draw (xj) node[below right] {\scriptsize{$\hat{p}_j(t)$}};
            \def\dis{0.5}
            \coordinate (p1) at (0,\dis);
            \coordinate (p2) at (\dis,\dis);
            \coordinate (p3) at (\dis,0);
            \draw [] (p1)--(p2)--(p3);
            \coordinate (p4) at (0,4-\dis);
            \coordinate (p5) at (\dis,4-\dis);
            \coordinate (p6) at (\dis,4);
            \draw [] (p4)--(p5)--(p6);
\end{tikzpicture}
  \caption{The geometric interpretation of protocol \eqref{eq_bearingNetworkLocalizeProtocal_element}.}
  \label{fig_protocolGeometricMeaning}
\end{figure}

Several remarks for protocol \eqref{eq_bearingNetworkLocalizeProtocal_element} are given below.
First, the protocol is distributed because the localization of $p_i$ only requires $\{g_{ij}\}_{j\in\N_i}$ and $\{\hat{p}_j\}_{j\in\N_i}$. In practical implementation, the bearings $\{g_{ij}\}_{j\in\N_i}$ can be measured by a bearing-only sensor such as a camera and the estimates $\{\hat{p}_j\}_{j\in\N_i}$ can be transmitted from the neighbors via wireless communication.
All the bearings must be measured in a global reference frame.
Second, the protocol has a clear geometric interpretation as shown in Figure~\ref{fig_protocolGeometricMeaning}.
The term $-P_{{g}_{ij}}(\hat{p}_i(t)-\hat{p}_j(t))$ is the orthogonal projection of $(\hat{p}_j(t)-\hat{p}_i(t))$ onto the orthogonal compliment of $g_{ij}$, and hence it acts to steer the estimate $\hat{p}_i(t)$ to align with the bearing measurement $g_{ij}$.
Third, protocol \eqref{eq_bearingNetworkLocalizeProtocal_element} can be viewed as an extension of the protocol proposed in \cite{ZhuGuangwei2014Automatica}, which is applicable to networks in the two-dimensional space.
Finally, those who are familiar with consensus problems might have noticed that protocol~\eqref{eq_bearingNetworkLocalizeProtocal} has a similar expression as the well-known consensus protocol \cite{OlfatiTAC2004}.
The difference is that in the consensus protocol, the weight for each edge is a positive scalar whereas in the localization protocol the weight for each edge is a positive semi-definite orthogonal projection matrix.

The convergence of the protocol is characterized as below.

\begin{theorem}\label{theorem_bearingNetworkProtocolConvergence}
The distributed protocol \eqref{eq_bearingNetworkLocalizeProtocal_element} can globally localize the network $\G(p)$ if and only if the network is localizable.
\end{theorem}
\begin{pf}
When $\B_{ff}$ is nonsingular (i.e., the network is localizable), the matrix $-\B_{ff}$ is Hurwitz.
As a result, the linear time-invariant system \eqref{eq_bearingNetworkLocalizeProtocal} is stable and the state converges to the steady state value $-\B_{ff}^{-1}\B_{fa}p_a$ which equals to the real follower location $p_f$ according to Lemma~\ref{lemma_pa_pf_relation}. When $\B_{ff}$ is singular (i.e., the network is not localizable), the final estimate would depend on the initial estimate of the network location.
\qed\end{pf}

\subsection{Sensitivity Analysis}\label{section_sensitivityAnalysis}

Since the bearing measurements may be corrupted by errors in practice, it is meaningful to study the impact of constant measurement errors on the localization protocol \eqref{eq_bearingNetworkLocalizeProtocal_element}.
Denote the unit vector $\tilde{g}_{ij}\in\R^d$ as the measurement of $g_{ij}$.
In the presence of bearing measurement errors, the localization protocol \eqref{eq_bearingNetworkLocalizeProtocal} becomes
\begin{align}\label{eq_bearingNetworkLocalizeProtocal_errorCase}
\dot{\hat{p}}_f(t)=-\tilde{\B}_{ff}\hat{p}_f(t)-\tilde{\B}_{fa}p_a,
\end{align}
where $\tilde{\B}_{ff}$ and $\tilde{\B}_{fa}$ are obtained from $\B_{ff}$ and $\B_{fa}$ by replacing $g_{ij}$ with $\tilde{g}_{ij}$, respectively.
The matrix $\tilde{\B}_{ff}$ may not be symmetric since $\tilde{g}_{ij}\ne-\tilde{g}_{ji}$ in general.

We next analyze two problems regarding \eqref{eq_bearingNetworkLocalizeProtocal_errorCase}.
The first is when $\tilde{\B}_{ff}$ is positive stable (i.e., all its eigenvalues have positive real parts) such that \eqref{eq_bearingNetworkLocalizeProtocal_errorCase} is globally stable.
If $\tilde{\B}_{ff}$ is positive stable, the final estimate given by \eqref{eq_bearingNetworkLocalizeProtocal_errorCase} is
\begin{align}\label{eq_finalEstimateWithBearingError}
\hat{p}_f^*=-\tilde{\B}_{ff}^{-1}\tilde{\B}_{fa}p_a.
\end{align}
The second problem is how large the localization error $\|\hat{p}^*_f-p_f\|$ is.
To solve the two problems, define
$$\Delta \B_{ff}\triangleq\tilde{\B}_{ff}-\B_{ff}, \quad \Delta \B_{fa}\triangleq\tilde{\B}_{fa}-\B_{fa},$$ as the perturbations of $\B_{ff}$ and $\B_{fa}$ caused by the bearing measurement errors.
Let $\theta_{ij}\in[0,\pi]$ be the angle between $\tilde{g}_{ij}$ and $g_{ij}$; that is $g_{ij}^\T\tilde{g}_{ij}=\cos\theta_{ij}$.
The angle $\theta_{ij}$ represents the inconsistency between $\tilde{g}_{ij}$ and $g_{ij}$.
This representation is valid for arbitrary dimensions.
Note $\theta_{ij}\ne\theta_{ji}$ in general.
Define the total bearing measurement error for the followers as
\begin{align*}
\epsilon\triangleq2\sum_{i\in\V_f}\sum_{j\in\N_i}\sin\theta_{ij}.
\end{align*}
We next give lemmas to characterize the relationship between $\epsilon$ and $\Delta \B_{ff}, \Delta\B_{fa}$.

\begin{lemma}\label{lemma_Px-Py_norm}
Denote by $\theta\in[0,\pi]$ the angle between any two nonzero vectors $x, y\in\R^d$ (i.e., $x^\T  y=\|x\|\|y\|\cos\theta$).  Then $\|P_x-P_y\|=\sin\theta.$
\end{lemma}
\begin{pf}
See Appendix~\ref{appendix_Proof of lemma_Px-Py_norm}.
\qed\end{pf}
\begin{lemma}\label{lemma_DeltaLffDeltaLfa_bearingError}
For a network $\G(p)$ with arbitrary bearing measurements $\{\tilde{g}_{ij}\}_{(i,j)\in\E}$, it always holds that $\|\Delta\B_{ff}\|\le \epsilon$ and $\|\Delta\B_{fa}\|\le \epsilon/2$.
\end{lemma}
\begin{pf}
Denote $\Delta P_{g_{ij}}\triangleq P_{\tilde{g}_{ij}}-P_{g_{ij}}, \forall (i,j)\in\E$.
It then follows from Lemma~\ref{lemma_Px-Py_norm} that $\|\Delta P_{g_{ij}}\|=\sin\theta_{ij}$.
Note $[\Delta \B_{ff}]_{ii}=\sum_{j\in\N_i}\Delta P_{g_{ij}}$ for $i\in\V_f$; $[\Delta \B_{ff}]_{ij}=-\Delta P_{g_{ij}}$ for $i\in\V_f$ and $j\in\N_i\cap\V_f$; and $[\Delta \B_{ff}]_{ij}=0$ otherwise.
Then we have
$
\|\Delta\B_{ff}\|
\le \sum_{i\in\V_f}\sum_{j\in\N_i\cap\V_f}\|\Delta P_{g_{ij}}\|+\sum_{i\in\V_f}\left\|\sum_{j\in\N_i}\Delta P_{g_{ij}}\right\|
\le \sum_{i\in\V_f}\sum_{j\in\N_i}\|\Delta P_{g_{ij}}\|+\sum_{i\in\V_f}\sum_{j\in\N_i} \|\Delta P_{g_{ij}}\|
\le 2\sum_{i\in\V_f}\sum_{j\in\N_i}\|\Delta P_{g_{ij}}\|
=2\sum_{i\in\V_f}\sum_{j\in\N_i}\sin\theta_{ij}=\epsilon.
$
Similarly, we have
$
\|\Delta\B_{fa}\|
\le \sum_{i\in\V_f}\sum_{j\in\N_i\cap\V_a}\|\Delta P_{g_{ij}}\|
\le \sum_{i\in\V_f}\sum_{j\in\N_i}\|\Delta P_{g_{ij}}\|
=\sum_{i\in\V_f}\sum_{j\in\N_i}\sin\theta_{ij}=\epsilon/2.
$
\qed\end{pf}

We now give a upper bound for the total bearing error $\epsilon$ to ensure the positive stability of $\tilde{\B}_{ff}$.

\begin{theorem}\label{theorem_sufficientConditionForMffNonsingular}
Given a localizable network with $\B_{ff}$ nonsingular, the matrix $\tilde{\B}_{ff}$ is positive stable if the total bearing error $\epsilon$ satisfies
\begin{align}\label{eq_sufficentConditionForMffNonsingular}
\epsilon<\lambda_{\min}(\B_{ff}),
\end{align}
where $\lambda_{\min}(\B_{ff})$ is the minimum eigenvalue of $\B_{ff}$.
\end{theorem}

\begin{pf}
Since $\|\Delta\B_{ff}\|<\epsilon$ by Lemma~\ref{lemma_DeltaLffDeltaLfa_bearingError},
if \eqref{eq_sufficentConditionForMffNonsingular} holds, we have $\|\Delta\B_{ff}\|<\lambda_{\min}(\B_{ff})=1/\|\B_{ff}^{-1}\|$,
which further implies $\|\B_{ff}^{-1}\Delta \B_{ff}\|\le\|\B_{ff}^{-1}\|\|\Delta \B_{ff}\|<1$.
Thus the spectral radius $\rho(\B_{ff}^{-1}\Delta\B_{ff})<1$ and hence the matrix $(I+\B_{ff}^{-1}\Delta\B_{ff})$ is nonsingular.
As a result, $\tilde{\B}_{ff}=\B_{ff}+\Delta\B_{ff}=\B_{ff}(I+\B_{ff}^{-1}\Delta\B_{ff})$ is nonsingular.
Since $\tilde{\B}_{ff}$ is obtained by perturbing $\B_{ff}$ and $\B_{ff}$ is positive stable, the nonsingularity of $\tilde{\B}_{ff}$ implies the positive stability.
\qed\end{pf}

Theorem~\ref{theorem_sufficientConditionForMffNonsingular} suggests that a large $\lambda_{\min}(\B_{ff})$ would give the network a large tolerance to bearing measurement errors.

We now study the localization error $\|\hat{p}_f^*-p_f\|$.
An intuitive conclusion that can be immediately drawn from \eqref{eq_finalEstimateWithBearingError} and matrix perturbation theory is that the localization error would be sufficiently small when the bearing measurement errors are sufficiently small.
We next give a specific upper bound on the localization error.

\begin{theorem}\label{theorem_upperBoundofEstimateError}
The estimate $\hat{p}^*_f=-\tilde{\B}_{ff}^{-1}\tilde{\B}_{fa}p_a$ given in \eqref{eq_finalEstimateWithBearingError} satisfies
$
\|\hat{p}_f^*-p_f\|
\le \frac{\epsilon}{\lambda_{\min}(\B_{ff})-\epsilon}\left(\frac{1}{2}\|p_a\|+\|p_f\|\right)$.
\end{theorem}
\begin{pf}
See Appendix~\ref{appendix_proof_theorem_upperBoundofEstimateError}.
\qed\end{pf}

\begin{figure*}
  \centering
  \subfloat[Initial estimate]{\includegraphics[width=0.2\linewidth]{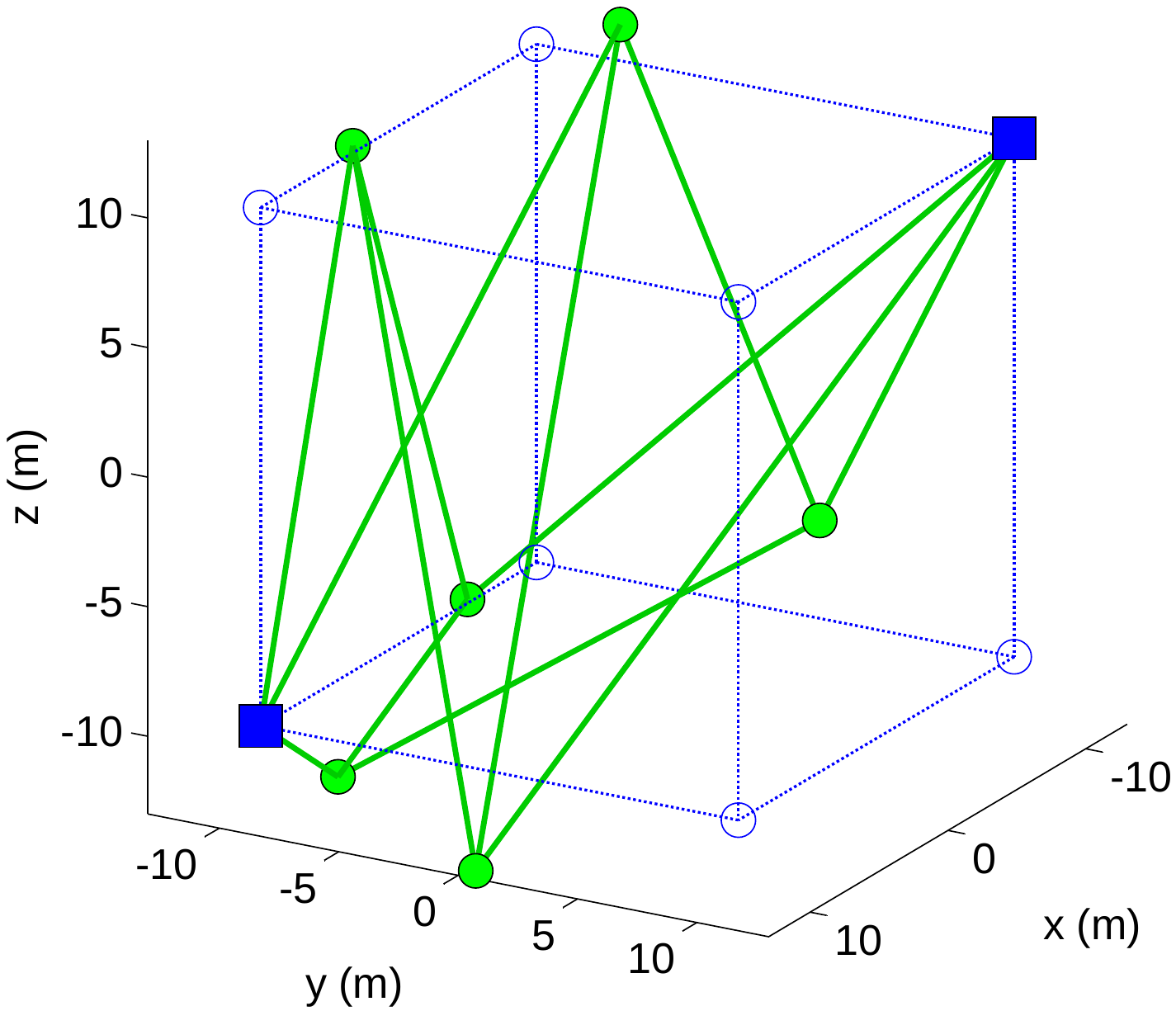}}
  \subfloat[$\sum_{i\in\V_f}\|\hat{p}_i(t)-p_i\|$]{\includegraphics[width=0.2\linewidth]{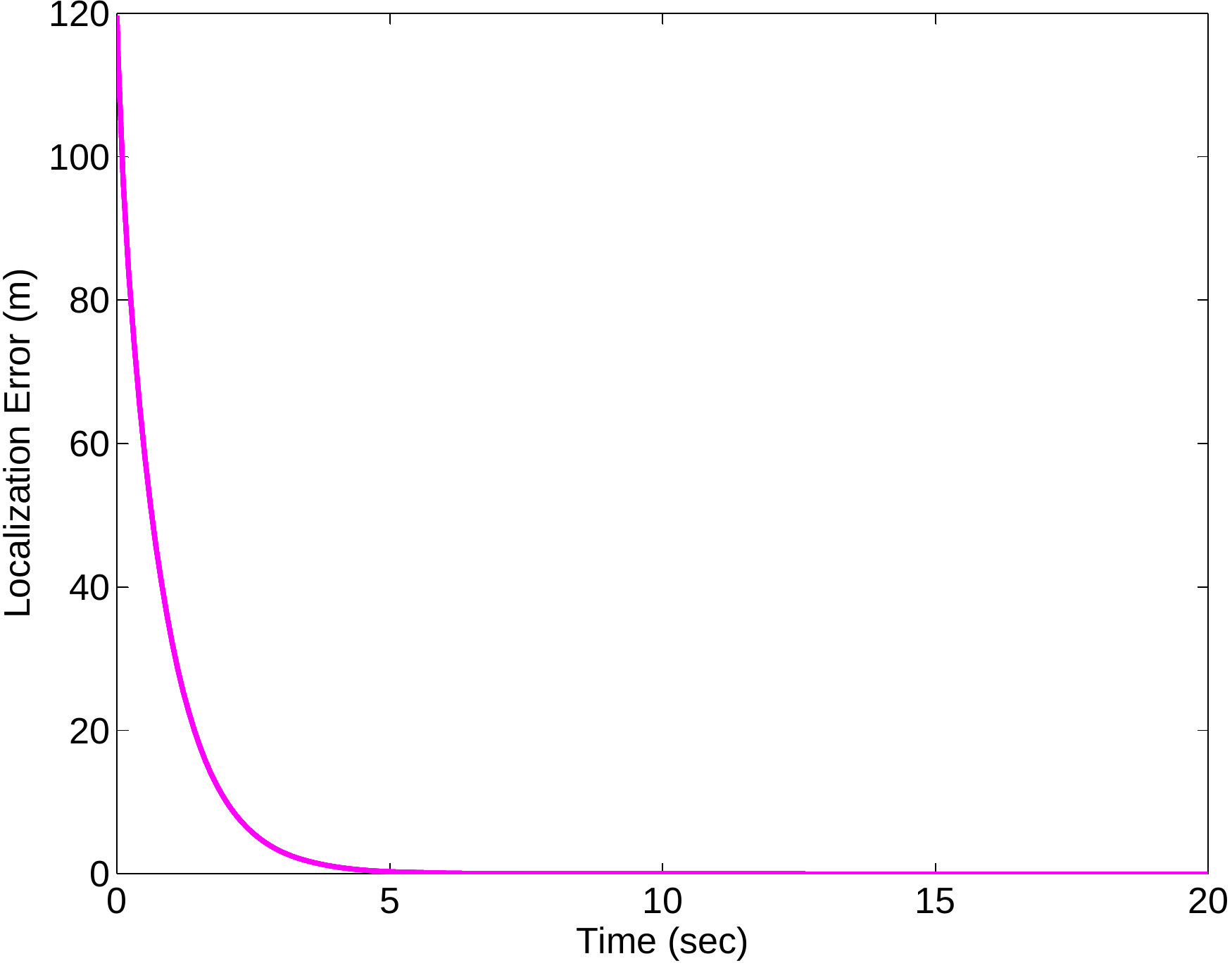}}
  \subfloat[Final estimate]{\includegraphics[width=0.2\linewidth]{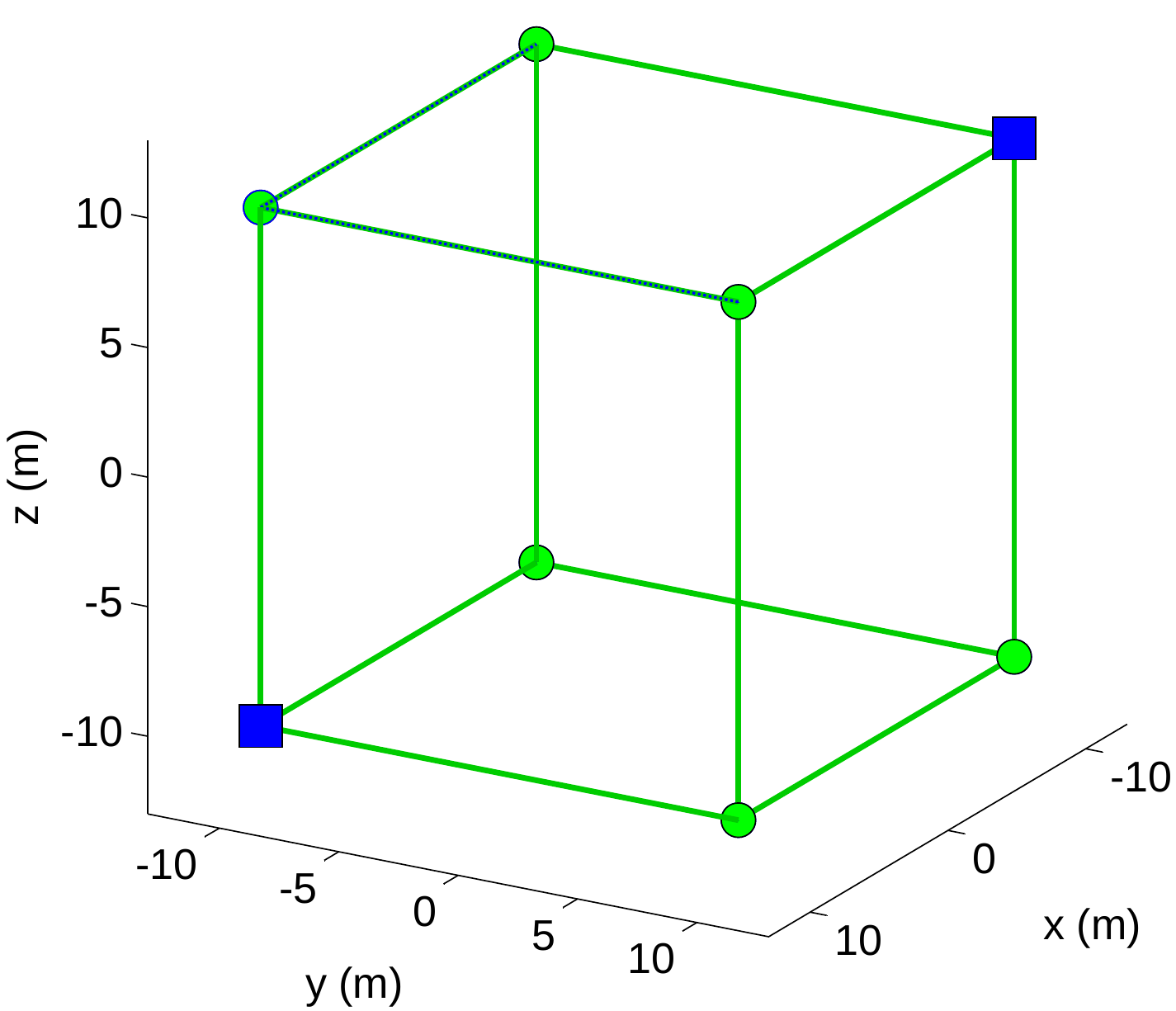}}
  \subfloat[$\sum_{i\in\V_f}\|\hat{p}_i(t)-p_i\|$]{\includegraphics[width=0.2\linewidth]{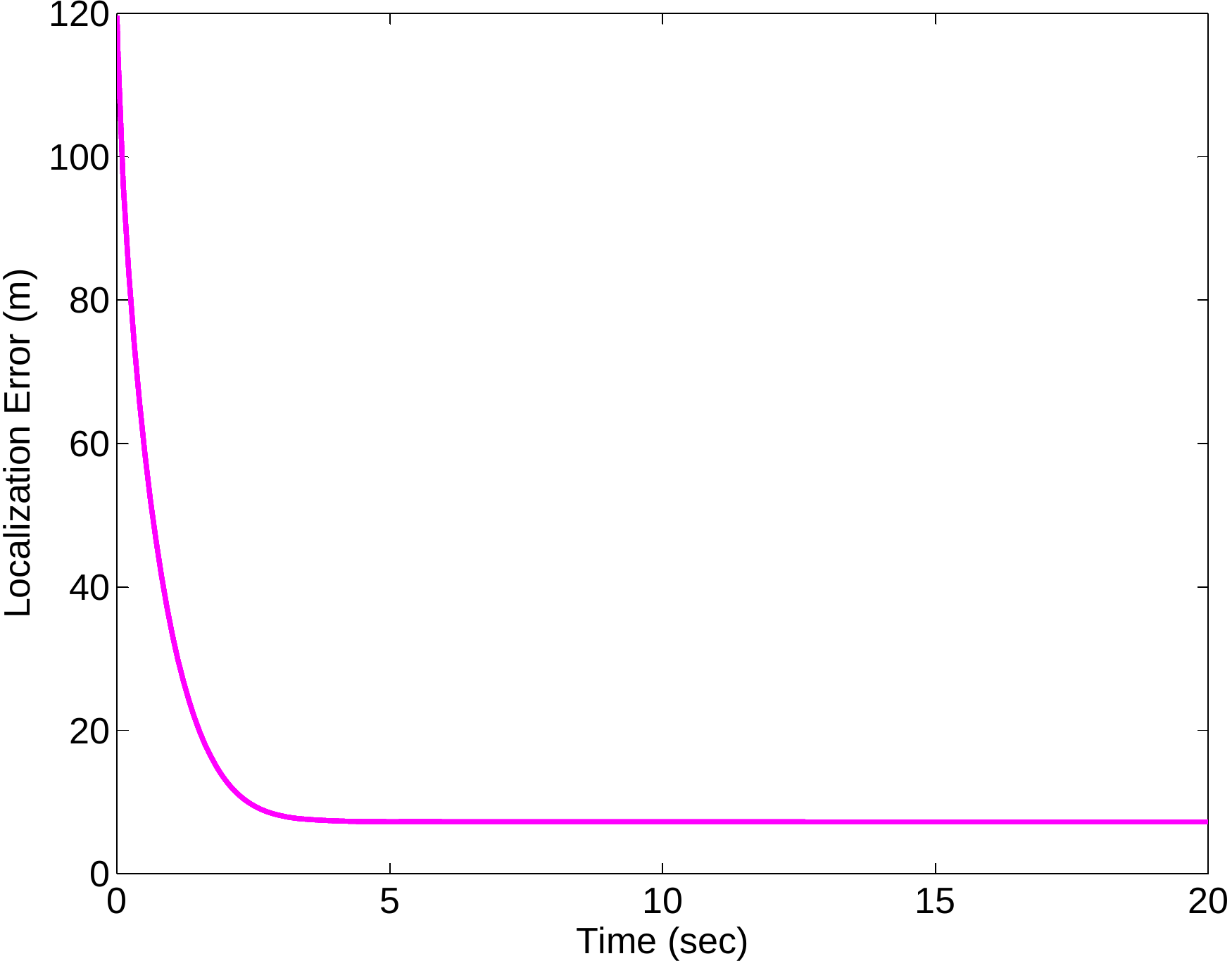}}
  \subfloat[Final estimate]{\includegraphics[width=0.2\linewidth]{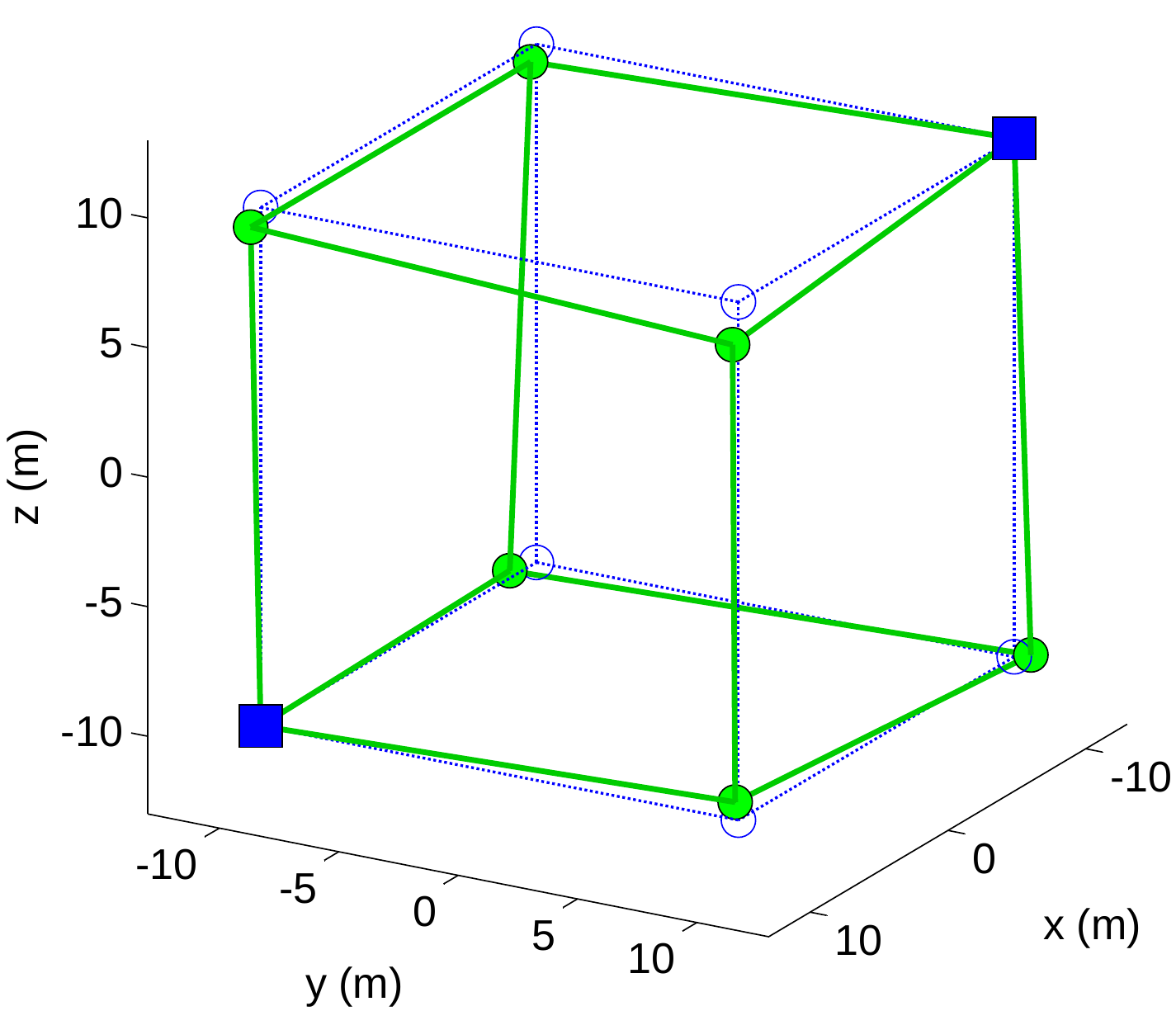}}
  \caption{Simulation examples for the localization protocol~\eqref{eq_bearingNetworkLocalizeProtocal_element}. The bearing measurements are \emph{accurate} for the example in (b)-(c), and \emph{inaccurate} for the one in (d)-(e).
  The blue squares represent the anchors. The blue hollow dots and the green solid dots represent the true and estimated locations of the followers, respectively.
  }
  \label{fig_sim_bearingLocalization}
\end{figure*}

In the last, we briefly discuss the impact of measurement errors in the anchors' locations.
Suppose the bearing measurements are accurate in this case.
Then the final estimate given by protocol \eqref{eq_bearingNetworkLocalizeProtocal_element} becomes $\hat{p}_f^*=-\B_{ff}^{-1}\B_{fa}(p_a+\Delta p_a)$, where $\Delta p_a\in\R^{dn_a}$ denotes the anchor location error.
Then the localization error is given by $\Delta \hat{p}_f\triangleq\hat{p}_f^*-p_f=-\B_{ff}^{-1}\B_{fa}\Delta p_a$, which indicates that the anchor location errors prorogate to the final localization error via a linear transformation.
It is straightforward to show that a translational or scaling error in the anchor measurements would cause the same translational or scaling error in the localization of followers.

\subsection{Simulation Examples}

Two simulation examples are shown in Figure~\ref{fig_sim_bearingLocalization} to demonstrate the localization protocol
\eqref{eq_bearingNetworkLocalizeProtocal_element}.
The network to be localized is a three-dimensional cubic network, which contains eight nodes and two of them are anchors and the other six are followers.
The initial estimate, which is randomly generated, is given in Figure~\ref{fig_sim_bearingLocalization}(a).
For the first example in Figure~\ref{fig_sim_bearingLocalization}(b)-(c), the bearing measurements are accurate and it can be seen that the estimate of the network location converges to the true value.
For the second example in Figure~\ref{fig_sim_bearingLocalization}(d)-(e), the bearing measurements are inaccurate.
Specifically, the total bearing error is $\epsilon=2.77$ and the final localization error equals $7.25$~m.
By comparing the two examples, it can be seen that when the bearings have measurement errors, the finally localized network would have localization errors. However, the final localized network can still be sufficiently close to the true network if the bearing errors are sufficiently small.
In addition, for the second example, we have $\lambda_{\min}=0.59<\epsilon$.
Although the condition in Theorem~\ref{theorem_sufficientConditionForMffNonsingular} is not satisfied, the matrix $\tilde{\B}_{ff}$ is still positive stable which indicates that the condition in Theorem~\ref{theorem_sufficientConditionForMffNonsingular} may be conservative.

\section{Conclusions}\label{section_conclusion}

This paper studied the problem of bearing-based network localization in arbitrary dimensions.
The first main contribution of this paper is to propose a variety of necessary and/or sufficient conditions for network localizability.
The second main contribution is to propose and analyze a linear localization protocol.
The results presented in this paper not only can be applied to solve the problem of sensor network localization but also provide a theoretical foundation for bearing-based multi-agent formation control \cite{zhao2015ECC,zhao2015MSC,zhao2015Maneuver,zhao2015CDC}.

In this paper, we assumed that the underlying graph is undirected.
As we have explained, the localizability analysis is independent to whether or not the sensing graph is undirected because any directed graph can be converted to an undirected one without affecting the localizability analysis.
But the convergence analysis of the proposed localization protocol relies on the assumption of undirected graphs.
As observed in \cite{zhao2015CDC}, a new notion termed bearing persistence emerges and makes the problem more complicated to analyze in the directed case.
Distributed localization with directed interaction topologies is therefore a direction for future work.

{\small
\section*{Acknowledgements}
The work presented here has been supported by the Israel Science Foundation (grant no. 1490/13).

\appendix

\section{Preliminaries to Bearing Rigidity Theory}\label{appendix_preliminaryBearingRigidity}

For a network $\G(p)$, consider an oriented graph and express the edge vector and the bearing for the $k$th directed edge in the oriented graph, respectively, as $e_{k}$ and $g_{k}\triangleq {e_{k}}/{\|e_{k}\|}$ for $k\in\{1,\dots,m\}$.
Define the \emph{bearing function} $F_B: \R^{dn}\rightarrow\R^{dm}$ as
$
F_B(p)\triangleq [g_1^\T ,\dots, g_m^\T]^\T.
$
The \emph{bearing rigidity matrix} is defined as the Jacobian of the bearing function,
$R_B(p) \triangleq \partial F_B(p)/\partial p\in\R^{dm\times dn}.$
Two important properties of the bearing rigidity matrix are given as below.

\begin{lemma}[\cite{zhao2014TACBearing}]\label{lemma_bearingRigidityMatrixRank}
For any network $\G(p)$, the bearing rigidity matrix satisfies
$R_B= \mydiag\left({P_{g_k}}/{\|e_k\|}\right)\bar{H}$, $\rank(R_B)\le dn-d-1$ and $\myspan{\one\otimes I_d, p}\subseteq \Null(R_B)$.
\end{lemma}

Let $\delta p$ be a variation of $p$.
If $R_B(p)\delta p=0$, then $\delta p$ is called an \emph{infinitesimal bearing motion} of $\G(p)$.
A network always has two kinds of \emph{trivial} infinitesimal bearing motions: translation and scaling of the entire network.

\begin{definition}[{Infinitesimal Bearing Rigidity}]\label{definition_infinitesimalParallelRigid}
    A network is \emph{infinitesimally bearing rigid} if all the infinitesimal bearing motions are trivial.
\end{definition}

The necessary and sufficient conditions for infinitesimal bearing rigidity are summarized as below.

\begin{theorem}[\cite{zhao2014TACBearing}]\label{theorem_conditionInfiParaRigid}
    For any network $\G(p)$, the following statements are equivalent:
\begin{enumerate}[(a)]
    \item $\G(p)$ is {infinitesimally bearing rigid};
    \item $\G(p)$ can be uniquely determined up to a translation and a scaling factor by the inter-neighbor bearings;
    \item $\rank(R_B)=dn-d-1$;
    \item $\Null(R_B)=\myspan{\one\otimes I_d, p}$.
    \end{enumerate}
\end{theorem}

\section{Proof of Lemma~\ref{lemma_Px-Py_norm}}\label{appendix_Proof of lemma_Px-Py_norm}

\begin{pf}
Here we only prove the case of $d=3$.
Without loss of generality, assume $x$ and $y$ are two unit vectors satisfying $\|x\|=\|y\|=1$.
Then, we have $P_x=I_d-xx^\T $, $P_y=I_d-yy^\T $, and hence
$
\|P_x-P_y\|=\|xx^\T -yy^\T \|.
$
There always exists an orthogonal matrix $U\in\R^{3\times 3}$ such that the two vectors $x$ and $y$ can be orthogonally transformed to $Ux=[1,0,,0]^\T$ and $Uy=[\cos\theta,\sin\theta,0]^\T$.
Since the spectral norm is invariant to orthogonal matrices, we have
\begin{align*}
\|P_x-P_y\|
&=\|U(xx^\T -yy^\T )U^\T \| \\
&=\left\|\left[
           \begin{array}{cc}
             1 & 0 \\
             0 & 0 \\
           \end{array}
         \right]-
         \left[
           \begin{array}{cc}
             \cos^2\theta & \cos\theta\sin\theta \\
             \sin\theta\cos\theta & \sin^2\theta\\
           \end{array}
         \right]
\right\|\\
&=\sin\theta\|Q\|,
\end{align*}
where
$         Q=\left[
           \begin{array}{cc}
             \sin\theta & -\cos\theta \\
             -\cos\theta & -\sin\theta\\
           \end{array}
         \right].
$
It is easy to see $Q^\T Q=I_2$ and hence $Q$ is an orthogonal matrix.
Then, $\|P_x-P_y\|=\sin\theta\|Q\|=\sin\theta\|I\|=\sin\theta$.
\qed\end{pf}
\section{Proof of Theorem~\ref{theorem_upperBoundofEstimateError}}\label{appendix_proof_theorem_upperBoundofEstimateError}

\begin{pf}
Recall $p_f=-\B_{ff}^{-1}\B_{fa}p_a$.
Rewrite $\hat{p}_f^*$ as $\hat{p}_f^*=-(\B_{ff}+\Delta\B_{ff})^{-1}(\B_{fa}+\Delta\B_{fa})p_a$.
By \cite[Eq.~(25)]{Henderson1981SumMatrix}, we have $(\B_{ff}+\Delta\B_{ff})^{-1}=\B_{ff}^{-1}-\B_{ff}^{-1}\Delta\B_{ff}(I+\B_{ff}^{-1}\Delta\B_{ff})^{-1}\B_{ff}^{-1},$
substituting which into $\hat{p}_f^*$ gives
$
\hat{p}_f^*
=-\B_{ff}^{-1}\B_{fa}p_a-\B_{ff}^{-1}\Delta\B_{fa}p_a +\B_{ff}^{-1}\Delta\B_{ff}(I+\B_{ff}^{-1}\Delta\B_{ff})^{-1}\B_{ff}^{-1}\Delta\B_{fa}p_a
+\B_{ff}^{-1}\Delta\B_{ff}(I+\B_{ff}^{-1}\Delta\B_{ff})^{-1}\B_{ff}^{-1}\B_{fa}p_a
=p_f-(I+\B_{ff}^{-1}\Delta\B_{ff})^{-1}\B_{ff}^{-1}\Delta\B_{fa}p_a
+\B_{ff}^{-1}\Delta\B_{ff}(I+\B_{ff}^{-1}\Delta\B_{ff})^{-1}p_f.
$
It follows that
\begin{align*}
\|\hat{p}_f^*-p_f\|
&\le \|(I+\B_{ff}^{-1}\Delta\B_{ff})^{-1}\B_{ff}^{-1}\Delta\B_{fa}p_a\| \\
&\qquad +\|\B_{ff}^{-1}\Delta\B_{ff}(I+\B_{ff}^{-1}\Delta\B_{ff})^{-1}p_f\| \\
&\le \|(I+\B_{ff}^{-1}\Delta\B_{ff})^{-1}\|\|\B_{ff}^{-1}\|\|\Delta\B_{fa}\|\|p_a\| \\
&\qquad +\|\B_{ff}^{-1}\|\|\Delta\B_{ff}\|\|(I+\B_{ff}^{-1}\Delta\B_{ff})^{-1}\|\|p_f\| \\
&=\|(I+\B_{ff}^{-1}\Delta\B_{ff})^{-1}\|\|\B_{ff}^{-1}\|  \left(\|\Delta\B_{fa}\|\|p_a\|\right. \\
&\qquad\qquad\qquad\qquad\qquad\qquad\qquad \left.+\|\Delta\B_{ff}\|\|p_f\|\right)
\end{align*}
Substituting $\|\Delta\B_{ff}\|\le \epsilon$ and $\|\Delta\B_{fa}\|\le\epsilon/2$ as shown in Lemma~\ref{lemma_DeltaLffDeltaLfa_bearingError}, and
$
\|(I+\B_{ff}^{-1}\Delta\B_{ff})^{-1}\|\le1/{(1-\|\B_{ff}^{-1}\|\|\Delta\B_{ff}\|)}
$
by \cite[Lemma~2.3.3]{BookMatrixComputation} into the above inequality gives
\begin{align*}
\|\hat{p}_f^*-p_f\|
&\le \frac{\|\B_{ff}^{-1}\|(\frac{1}{2}\|p_a\|+\|p_f\|)\epsilon}{1-\|\B_{ff}^{-1}\|\epsilon}.
\end{align*}
Substituting $\|\B_{ff}^{-1}\|=1/\lambda_{\min}(\B_{ff})$ completes the proof.
\qed\end{pf}

\bibliographystyle{ieeetr}      
\bibliography{myOwnPub,zsyReferenceAll}


\end{document}